\documentclass[12pt,a4paper]{article}
\usepackage{amsthm}
\begin{document}
{\Large Another method to solve the grasshopper problem\\(the International Mathematical Olympiad) } \\\\
\centerline{Yasushi Ieno (vei04156@nifty.com)}\\
\\\\
\centerline{Abstract}

The 6th problem of the 50th International Mathematical Olympiad (IMO), held in Germany, 2009, is called 'the grasshopper problem'. To this problem Kos[1] developed theory from unique viewpoints by reference of Noga Alon’s combinatorial Nullstellensatz.

We have tried to solve this problem by an original method inspired by a polynomial function that Kos defined in [1] , then examined for n=3, 4 and 5. For almost cases the claim of this problem follows, but there remains imperfection due to 'singularity'.\\

\noindent 0.Introduction\\

The 6th problem of the 50th International Mathematical Olympiad (IMO), held in Germany,
2009, was the following.\\

Let $a_1, a_2,\ldots, a_n$ be distinct positive integers and let M be a set of n$-$1 positive integers not containing s=$a_1$+$a_2$$+$$\cdots$$+$$a_n$. A grasshopper is to jump along
the real axis, starting at the point 0 and making n jumps to the right with lengths
$a_1, a_2,\ldots, a_n$ in some order. Prove that the order can be chosen in such a way that
the grasshopper never lands on any point in M.\\

According to [1], Kos says that up to now, all known solutions to this problem, so called 'the grasshopper problem', are elementary and inductive, for example, by drawing a real axis on paper. In fact a solution of ours below is one of its examples. 

Then in [1], Kos tried to apply Noga Alon's combinatorial Nullstellensatz [2], which is effective but not perfect to solve the grasshopper problem, as a result he could not solve the problem with his intentional method. 

So we try to present a way to solve the problem and prove it by reference of [1], even if partially. \\
\\\\
\noindent 1. Alon's combinatorial Nullstellensatz\\

Now we introduce an interesting tool which may help our investigation.\\
\\
Lemma 1 (Nonvanishing combinatorial Nullstellensatz). 

Let $S_1$,\ldots,$S_n$ be
nonempty subsets of a field F, and let $t_1$,\ldots,$t_n$ be nonnegative integers such that $t_i$$<$$|S_i|$
for i=$1$,$2$,$\ldots$,n. Let P($x_1$,$\ldots$,$x_n$) be a polynomial over F with total degree $t_1$+$\cdots$+$t_n$,
and suppose that the coefficient of ${x_1}^{t_1}{x_2}^{t_2}$$\cdots$${x_n}^{t_n}$
in P($x_1$,\ldots,$x_n$) is nonzero. Then there exist
elements $s_1$$\in$$S_1$,\ldots,$s_n$$\in$$S_n$ for which P($s_1$,\ldots,$s_n$)$\neq$0.\\

Also we present a polynomial function {\it f{\rm (}$x_1$,$x_2$,$\ldots,$$x_n$}) by reference of  [1] as follows.
\begin{eqnarray}
&& f(x_1,x_2,\ldots,x_n):=(x_1{\rm-}m_1{\rm)(}x_1{\rm-}m_2)\ldots(x_1{\rm-}m_{n-1})(x_1{\rm+}x_2{\rm-}m_1)\ldots \nonumber \\
&& \ \ \ \ \ \ (x_1{\rm+}x_2{\rm-}m_2)\cdots(x_1{\rm+}x_2{\rm-}m_{n-1})(x_1{\rm+}\cdots{\rm+}x_{n-1}{\rm-}m_1)\cdots \nonumber  \\
&& \ \ \ \ \ \ (x_1{\rm+}\cdots{\rm+}x_{n-1}{\rm-}m_2) \cdots(x_1{\rm+}\cdots{\rm+}x_{n-1}{\rm-}m_{n-1}) \nonumber \\
&& \ \ \ \ \ \ =\prod_{l=1}^{n-1} {\prod_{i=1}^{n-1}{((x_1{\rm+}x_2{\rm+}\cdots{\rm+}x_l){\rm-}m_i)}} 
\end{eqnarray}

On the grasshopper problem now if we fix the jumping order as $a_1$,$a_2,$$\ldots$,$\\
$$a_n$, then a grasshopper succeeds in its jumping without blocked if and only if f($a_1$,$...$,$a_n$)$\neq$0, then the degree of f($a_1$,$...$,$a_n$) is $({\rm n}$$-$$1)^2$. And $x_1^{\rm n-1}x_2^{\rm n-1}\cdots$\\
$x_{n-1}^{\rm n-1}$ is a monomial the total degree of which is $({\rm n}$$-$$1)^2$, and the coefficient of which is 1.

Now we define n sets $S_1$,$S_2$,$\ldots$,$S_n$ such that $S_1$$=$$S_2$$=$$\cdots$$=$$S_n$$=$$\{$$a_1$,$a_2$,$\ldots,$\\
$a_n$$\}$, then the number of elements of these n sets are $|S_1|$$=$$|S_2|$$=$$\cdots$$=$$|S_n|$$=$${\rm n}$$\\
>$n$-$1, so we can adopt Lemma 1 to this polynomial function (1).

But there remains imperfection because the elements $a_1$,$a_2$,$...$,$a_n$ considered in Lemma 1 are not necessarily distinct, that is to say, a pair of ($a_1$,\ldots,$a_n$) may be the same number. 

If we multiple f($x_1$,$\ldots$,$x_n$) by the so-called Vandermonde polynomial (see, for example, [3, pp. 346–347]), a new polynomial is created as follows. \\
\begin{eqnarray}
&& \prod_{1\le k<j\le {\rm n-1}}^{ }{(x_k{-}x_j)} \prod_{l=1}^{\rm n-1} {\prod_{i=1}^{\rm n-1}{((x_1{+}x_2{+}\cdots{+}x_l){-}m_i)}}
\end{eqnarray}

The elements $a_1$,$a_2$,$\ldots$,$a_n$ are required to be distinct if the new polynomial is nonzero when $x_i$=$a_i$ for any i such that 1$\leq$i$\leq$n. But any monomial of (2) the total degree of which is equal to the degree of (2), $({\rm n}{\rm -1})^2$+${}_{{\rm n}{\rm -1}} {\rm C}_2$, has a factor the exponent of which is over {\rm n}$-$1. Thus Lemma 1 can not be applied.\\
\\
\noindent 2. Attempts to use new polynomials by permutations\\

We could not apply Lemma 1 to f($x_1$,$x_2$,\ldots,$x_n$) if $a_1$,$a_2$,$\ldots$,$a_n$ are distinct.

We want to find out an effective polynomial function, on the condition that the total degree is kept, if possible. 

Let ${\rm Sym}({\rm n})$ be a symmetric group of degree n.
By a permutation ${\pi}$$\in$${\rm Sym}({\rm n})$, we get
\begin{eqnarray}
&& f(x_{\pi(1)},x_{\pi(2)},\ldots,x_{\pi(n)}) \nonumber \\
&&\ \ \ \ ={\prod_{l=1}^{\rm n-1} {\prod_{i=1}^{{\rm n}{\rm-1}}{((x_{\pi(1)}{+}x_{\pi(2)}{+}\cdots{+}x_{\pi(l)}){-}m_i)}}}
\end{eqnarray}

There are totally $({\rm n}{\rm -1})^2$ factors in (3). 

And the total number of cases by possible permutations is n!.

Then we multiple each (3) by the signature of each permutation, that is +1 or $-$1, and make their summation as follows.
\begin{eqnarray}
&& {\sum_{\pi \in {\rm Sym}(n)}^{ }} {\rm sgn}(\pi) f(x_{\pi(1)},x_{\pi(2)},\ldots,x_{\pi(n)}) \nonumber \\
&& \ \ \ \ ={\sum_{\pi \in {\rm Sym}(n)}^{ } {\rm sgn}(\pi) \prod_{l=1}^{n-1} {\prod_{i=1}^{n-1} {((x_{\pi(1)}{+}x_{\pi(2)}+\cdots{+}x_{\pi(l)}){-}m_i)}}}
\end{eqnarray}
\\

In (4) $x_i$ and $x_j$ is anti-symmetric if i is not equal to j, so it may be a multiple of the above-mentioned Vandermonde polynomial.
\\\\
3. Real example for this case\\
\\
3-1. the case n=3\\

Unfortunately Alon's combinatorial Nullstellensatz can't be applied now, because by simple computations we can see that nothing but unsuitable 4-degree monomials like $x_1^3x_2$, $x_1^3x_3$ exist. In this case $|S_1|$ must be larger than 3, applying Lemma 1 is impossible.

We compute (4) for n=3 by summing up 3!=6 polynomials as follows.
\begin{eqnarray}
&& {\sum_{\pi \in{\rm Sym(3)}}^{ } {\rm sgn}(\pi) f(x_{\pi(1)},x_{\pi(2)})} \nonumber \\
&& \ \ \ \ ={\sum_{\pi \in{\rm Sym(3)}}^{ } {\rm sgn}(\pi) \prod_{l=1}^{2} {\prod_{i=1}^{2} {((x_{\pi(1)}{+}x_{\pi(2)}){-}m_i)}}}
\end{eqnarray}

The computation of (5) is the following.
\begin{eqnarray}
&& (5)\ {=}\ f(x_{1},x_{2},x_{3}){-}f(x_{1},x_{3},x_{2}){-}f(x_{2},x_{1},x_{3})  \nonumber \\
&& \ \ \ \ \ \ {+}f(x_{2},x_{3},x_{1}){+}f(x_{3},x_{1},x_{2}){-}f(x_{3},x_{2},x_{1})  \nonumber \\
&& \ \ \ \ \ \ {=}{({x_1}{-}{m_1})({x_1}{-}{m_2})({x_1}{+}{x_2}{-}{m_1})({x_1}{+}{x_2}{-}{m_2})} \nonumber \\
&& \ \ \ \ \ \ {-}{({x_1}{-}{m_1})({x_1}{-}{m_2})({x_1}{+}{x_3}{-}{m_1})({x_1}{+}{x_3}{-}{m_2})} \nonumber \\
&& \ \ \ \ \ \ {-}{({x_2}{-}{m_1})({x_2}{-}{m_2})({x_2}{+}{x_1}{-}{m_1})({x_2}{+}{x_1}{-}{m_2})} \nonumber \\
&& \ \ \ \ \ \ {+}{({x_2}{-}{m_1})({x_2}{-}{m_2})({x_2}{+}{x_3}{-}{m_1})({x_2}{+}{x_3}{-}{m_2})} \nonumber \\
&& \ \ \ \ \ \ {+}{({x_3}{-}{m_1})({x_3}{-}{m_2})({x_3}{+}{x_1}{-}{m_1})({x_3}{+}{x_1}{-}{m_2})} \nonumber \\
&& \ \ \ \ \ \ {-}{({x_3}{-}{m_1})({x_3}{-}{m_2})({x_3}{+}{x_2}{-}{m_1})({x_3}{+}{x_2}{-}{m_2})} \nonumber \\
&& \ \ \ \ \ \ {=}{({x_1}{-}{x_2})({x_1}{-}{x_3})({x_2}{-}{x_3})(({x_1}{+}{x_2}{+}{x_3}){-}({m_1}{+}{m_2}))}
\end{eqnarray}

We present other computations. 3 pairs of the above 6 polynomials appear by turns.
\begin{eqnarray}
&& \ \ f(x_{1},x_{2},x_{3}){-}f(x_{1},x_{3},x_{2}) \nonumber \\
&& \ \ \ \ {=}{({x_1}{-}{m_1})({x_1}{-}{m_2})({x_1}{+}{x_2}{-}{m_1})({x_1}{+}{x_2}{-}{m_2})} \nonumber \\
&& \ \ \ \ {-}{({x_1}{-}{m_1})({x_1}{-}{m_2})({x_1}{+}{x_3}{-}{m_1})({x_1}{+}{x_3}{-}{m_2})} \nonumber \\
&& \ \ \ \ {=}{({x_1}{-}{m_1})({x_1}{-}{m_2})((2{x_1}{+}{x_2}{+}{x_3}){-}({m_1}{+}{m_2}))} \\
&& \ \ f(x_{2},x_{1},x_{3}){-}f(x_{2},x_{3},x_{1}) \nonumber \\
&& \ \ \ \ {=}{({x_2}{-}{m_1})({x_2}{-}{m_2})({x_2}{+}{x_1}{-}{m_1})({x_2}{+}{x_1}{-}{m_2})} \nonumber \\
&& \ \ \ \ {-}{({x_2}{-}{m_1})({x_2}{-}{m_2})({x_2}{+}{x_3}{-}{m_1})({x_2}{+}{x_3}{-}{m_2})} \nonumber \\
&& \ \ \ \ {=}{({x_2}{-}{m_1})({x_2}{-}{m_2})(({x_1}{+}2{x_2}{+}{x_3}){-}({m_1}{+}{m_2}))} \\
&& \ \ f(x_{3},x_{1},x_{2}){-}f(x_{3},x_{2},x_{1}) \nonumber \\
&& \ \ \ \ {=}{({x_3}{-}{m_1})({x_3}{-}{m_2})({x_3}{+}{x_1}{-}{m_1})({x_3}{+}{x_1}{-}{m_2})} \nonumber \\
&& \ \ \ \ {-}{({x_3}{-}{m_1})({x_3}{-}{m_2})({x_3}{+}{x_2}{-}{m_1})({x_3}{+}{x_2}{-}{m_2})} \nonumber \\
&& \ \ \ \ {=}{({x_3}{-}{m_1})({x_3}{-}{m_2})(({x_1}{+}{x_2}{+}2{x_3}){-}({m_1}{+}{m_2}))} 
\end{eqnarray}

\noindent Theorem 1.
\begin{eqnarray}
&& \ \ \ {\rm Let \ } a_1, a_2, a_3\ {\rm be\ distinct\ positive\ integers,\ and\ }m_1, m_2\ {\rm be\ distinct} \nonumber \\
&& {\rm positive\ integers,\ then\ there\ exists\ } {\pi \in{\rm Sym(3)}}\ {\rm that\ holds} \nonumber \\
&& \ \ \ f(a_{\pi(1)},a_{\pi(2)},a_{\pi(3)})= \nonumber \\
&& \ \ \ \ (a_{\pi(1)}{-}m_1)(a_{\pi(1)}{-}m_2)(a_{\pi(1)}{+}a_{\pi(2)}{-}m_1)(a_{\pi(1)}{+}a_{\pi(2)}{-}m_2) \nonumber \\
&& \ \ \ \ \ \ \ \ \ \neq 0.
\end{eqnarray}
\begin{proof}
\quad\par
If\ \ \ $f(a_{\pi(1)},a_{\pi(2)})$=$(a_{\pi(1)}{-}m_1)(a_{\pi(1)}{-}m_2)(a_{\pi(1)}{+}a_{\pi(2)}{-}m_1)$\\
$\times(a_{\pi(1)}{+}a_{\pi(2)}{-}m_2)$=0 for any ${\pi{\in}{\rm Sym(3)}}$,\\
then four equations hold as below by (6),(7),(8) and(9). 
\begin{eqnarray}
&& {({a_1}{-}{a_2})({a_1}{-}{a_3})({a_2}{-}{a_3})(({a_1}{+}{a_2}{+}{a_3}){-}({m_1}{+}{m_2}))}=0. \\
&& {({a_1}{-}{m_1})({a_1}{-}{m_2})((2{a_1}{+}{a_2}{+}{a_3}){-}({m_1}{+}{m_2}))}=0.\\
&& {({a_2}{-}{m_1})({a_2}{-}{m_2})(({a_1}{+}2{a_2}{+}{a_3}){-}({m_1}{+}{m_2}))}=0.\\
&& {({a_3}{-}{m_1})({a_3}{-}{m_2})(({a_1}{+}{a_2}{+}3{a_3}){-}({m_1}{+}{m_2}))}=0.
\end{eqnarray}

From (11), (${a_1}{+}{a_2}{+}{a_3}){-}({m_1}{+}{m_2}$)=0 follows, because $a_1, a_2, a_3$ are distinct. Then neither 2(${a_1}{+}{a_2}{+}{a_3}){-}({m_1}{+}{m_2}$) nor (${a_1}{+}2{a_2}{+}{a_3}){-}({m_1}{+}{m_2}$) nor (${a_1}{+}{a_2}{+}3{a_3}){-}({m_1}{+}{m_2}$) is equal to 0, so $({a_1}{-}{m_1})({a_1}{-}{m_2})$=0 \\
and $({a_2}{-}{m_1})({a_2}{-}{m_2})$=0 and $({a_3}{-}{m_1})({a_3}{-}{m_2})$=0 at (12), (13) and (14), which does not happen at the same time, this is because $a_1, a_2$ and $a_3$ are distinct and $m_1$ and $m_2$ are also distinct.  

It follows that the assumption above does not come true.

This completes the proof.

\end{proof}

If $f(a_1, a_2, a_3)\neq$0, then at least one of the above-mentioned six polynomials consisting of (6) is not 0. Therefore the claim of the grasshopper problem follows for n=3, that is to say, a grasshopper succeeds in jumping without landing on $m_1$ or $m_2$ by choosing one order $(a_{i1}, a_{i2}, a_{i3})$ out of six possible jumping orders, such that $f(a_{i1}, a_{i2}, a_{i3})$=$(a_{i1}{-}{m_1})({a_{i1}}{-}{m_2})({a_{i1}}{+}{a_{i2}}{-}{m_1})\\
({a_{i1}}{+}{a_{i2}}{-}{m_2})$$\neq$0. 

For the n=3's case of the grasshopper problem, $\{(a_1,a_2,a_3)|$$({a_1}{+}{a_2}{+}{a_3})$\\
${-}({m_1}{+}{m_2})=0$$\}$ is a 'singularity' set that may vanish the possibility of a grasshopper's safe jumping. But by comparing (6) with (7),(8) and (9), this possibility bas been easily denied. \\
\\
3-2.  the case n=4\\

We sum up 4!=24 polynomials which were made by permutation as follows.
\begin{eqnarray}
&& {\sum_{\pi \in{\rm Sym(4)}}^{ } {\rm sgn}(\pi) f(x_{\pi(1)},x_{\pi(2)},x_{\pi(3)})} \nonumber \\
&& \ \ \ \ ={\sum_{\pi \in{\rm Sym(4)}}^{ } {\rm sgn}(\pi) \prod_{l=1}^{3} {\prod_{i=1}^{3} {((x_{\pi(1)}{+}x_{\pi(2)}{+}x_{\pi(3)}){-}m_i)}}}
\end{eqnarray}

The degree is $3^2$=9 and the permutation number is 4!=24, so the computation of (15) is more complicated. We present the computing results for the case n=4, similarly as the case n=3, as below.
\begin{eqnarray}
&& (15)={(x_1{-}x_2)(x_1{-}x_3)(x_1{-}x_4)(x_2{-}x_3)(x_2{-}x_4)(x_3{-}x_4)} \nonumber \\
&& \ \times (3(x_1{+}x_2{+}x_3{+}x_4){-}2(m_1{+}m_2{+}m_3)) \nonumber \\
&& \ \times (6(x_1^2{+}x_2^2{+}x_3^2{+}x_4^2){+}8({m_1}{m_2}{+}{m_1}{m_3}{+}{m_1}{m_4}{+}{m_2}{m_3}{+}{m_2}{m_4}{+}{m_3}{m_4}) \nonumber \\
&& \ \ \ \ \ \ {-}7(m_1{+}m_2{+}m_3)(x_1{+}x_2{+}x_3{+}x_4) \nonumber \\
&& \ \ \ \ \ \ {+}(m_1^2{+}m_2^2{+}m_3^2{+}6{m_1}{m_2}{+}6{m_2}{m_3}{+}6{m_3}{m_1}))
\end{eqnarray}
\begin{eqnarray}
&& \ \ f(x_{1},x_{2},x_{3},x_{4}){-}f(x_{1},x_{2},x_{4},x_{3}) \nonumber \\
&& \ \ {=}({x_1}{-}{m_1})({x_1}{-}{m_2})({x_1}{-}{m_3})({x_1}{+}{x_2}{-}{m_1})({x_1}{+}{x_2}{-}{m_2})({x_1}{+}{x_2}{-}{m_3}) \nonumber \\
&& \ \ \times {({x_3}{-}{x_4})} \nonumber \\
&& \ \ \times({(3{x_1}^2{+}3{x_2}^2{+}{x_3}^2{+}{x_4}^2){+}(6{x_1}{x_2}{+}3{x_1}{x_3}{+}3{x_1}{x_4}{+}3{x_2}{x_3}{+}3{x_2}{x_4}{+}{x_3}{x_4})} \nonumber \\
&& \ \ \ \ \ \ {-}({m_1}{+}{m_2}{+}{m_3})(2{x_1}{+}2{x_2}{+}{x_3}{+}{x_4}){+}{m_1}{m_2}{+}{m_1}{m_3}{+}{m_2}{m_3})
\end{eqnarray}

Now generalizing (17), for $(x_{j1}$,$x_{j2}$,$x_{j3}$,$x_{j4})$, any permutation of $(x_{1}$,$x_{2}$,$x_{3}$,$x_{4})$, we obtain 
\begin{eqnarray}
&& \ \ f(x_{j1}{,}x_{j2}{,}x_{j3}{,}x_{j4}){-}f(x_{j1}{,}x_{j2}{,}x_{j4}{,}x_{j3}) \nonumber \\
&& \ \ {=}({x_{j1}}{-}{m_1})({x_{j1}}{-}{m_2})({x_{j1}}{-}{m_3})({x_{j1}}{+}{x_{j2}}{-}{m_1})({x_{j1}}{+}{x_{j2}}{-}{m_2})({x_{j1}}{+}{x_{j2}}{-}{m_3}) \nonumber \\
&& \ \ \times{(x_{j3}{-}x_{j4})} \nonumber \\
&& \ \ \times((3x_{j1}^2{+}3x_{j2}^2{+}x_{j3}^2{+}x_{j4}^2){+}(6x_{j1}x_{j2}{+}3x_{j1}x_{j3}{+}3x_{j1}x_{j4}{+}3x_{j2}x_{j3}{+}3x_{j2}x_{j4}{+}x_{j3}x_{j4}) \nonumber \\
&& \ \ \ \ {-}(m_1{+}m_2{+}m_3)(2x_{j1}{+}2x_{j2}{+}x_{j3}{+}x_{j4}){+}m_1m_2{+}m_1m_3{+}m_2m_3)
\end{eqnarray}

From (16), for the case n=4 of the grasshopper problem, we can obtain that
\begin{eqnarray}
&& \{(a_1,a_2,a_3,a_4)|\ \nonumber \\
&& \ \ (3(a_1{+}a_2{+}a_3{+}a_4){\rm -}2(m_1{+}m_2{+}m_3)) \nonumber \\
&& \ \ \times (6(a_1^2{+}a_2^2{+}a_3^2{+}a_4^2) \nonumber \\
&& \ \ \ \ {+}8({a_1}{a_2}{+}{a_1}{a_3}{+}{a_1}{a_4}{+}{a_2}{a_3}{+}{a_2}{a_4}{+}{a_3}{a_4}) \nonumber \\
&& \ \ \ \ {\rm -}7(m_1{+}m_2{+}m_3)(a_1{+}a_2{+}a_3{+}a_4)\nonumber \\
&& \ \ \ \ {+}(m_1^2{+}m_2^2{+}m_3^2{+}6{m_1}{m_2}{+}6{m_2}{m_3}{+}6{m_3}{m_1})){=}0\} \ \  
\end{eqnarray}
is a 'singularity' set that may eliminate the possibility of a grasshopper's safe jumping. 

Unlike the case n=3, the comparison of (18) and (19) does not lead to the solution of the grasshopper problem yet, for n=4.

For (17), when $({x_1}{,}{x_2}{,}{x_3}{,}{x_4})$=(1{,}4{,}2{,}3) and $({m_1}{,}{m_2}{,}{m_3})$=(2{,}3{,}10), then\\
$(3{x_1}^2{+}3{x_2}^2{+}{x_3}^2{+}{x_4}^2){+}(6{x_1}{x_2}{+}3{x_1}{x_3}{+}3{x_1}{x_4}{+}3{x_2}{x_3}{+}3{x_2}{x_4}{+}{x_3}{x_4})$\\
${-}({m_1}{+}{m_2}{+}{m_3})(2{x_1}{+}2{x_2}{+}{x_3}{+}{x_4}){+}{m_1}{m_2}{+}{m_1}{m_3}{+}{m_2}{m_3}$\\
=(3$\times$$1^2${+}3$\times$$4^2${+}$2^2${+}$3^2$)\\
{+}(6$\times$1$\times$4{+}3$\times$1$\times$2{+}3$\times$1$\times$3{+}3$\times$4$\times$2{+}3$\times$4$\times$3{+}2$\times$3)\\
${-}(2{+}3{+}10)(2$$\times$1{+}2$\times$4{+}2{+}3){+}(2$\times$3{+}2$\times$10{+}3$\times$10)\\
=64+105${-}$225+56=0

It follows that (17) is equal to 0 for this case, which does not require $({x_1}{-}{m_1})({x_1}{-}{m_2})({x_1}{-}{m_3})({x_1}{+}{x_2}{-}{m_1})({x_1}{+}{x_2}{-}{m_2})$$({x_1}{+}{x_2}{-}{m_3})$=0, \\
in fact, (1${-}$2)(1${-}$3)(1${-}$10)(1{+}4${-}$2)(1{+}4${-}$3)(1{+}4${-}$10)$\neq$0.

And the condition that $({x_1}{,}{x_2}{,}{x_3}{,}{x_4})$=(1{,}4{,}2{,}3) and $({m_1}{,}{m_2}{,}{m_3})$\\
=(2{,}3{,}10) fulfills (16)=0. In short, when $({m_1}{,}{m_2}{,}{m_3})$=(2{,}3{,}10), $({x_1}{,}{x_2}{,}{x_3}{,}{x_4})$\\
=(1{,}4{,}2{,}3) is an element of so-called the 'singularity' set. 

 $({x_1}{,}{x_2}{,}{x_3}{,}{x_4})$=(1{,}4{,}2{,}3) does not restrict the value of $({x_1}{-}{m_1})({x_1}{-}{m_2})({x_1}{-}{m_3})$\\
$\times({x_1}{+}{x_2}{-}{m_1})({x_1}{+}{x_2}{-}{m_2})({x_1}{+}{x_2}{-}{m_3})$, therefore we can not approach the proof of the case n=4 like Theorem 1.
\\\\
3-3.  the case n=5\\

We sum up 5!=120 polynomials which were made by permutation, as follows.
\begin{eqnarray}
&& {\sum_{\pi \in{\rm Sym(5)}}^{ } {\rm sgn}(\pi) f(x_{\pi(1)},x_{\pi(2)},x_{\pi(3)},x_{\pi(4)})} \nonumber \\
&& \ \ \ \ ={\sum_{\pi \in{\rm Sym(5)}}^{ } {\rm sgn}(\pi) \prod_{l=1}^{4} {\prod_{i=1}^{4} {((x_{\pi(1)}{+}x_{\pi(2)}{+}x_{\pi(3)}{+}x_{\pi(4)}){-}m_i)}}}
\end{eqnarray}

The computation of (20) is so complicated that we use the computational tool SageMath[5].

We only show the result in the appendix. 

Unlike the cases n=3 and 4, the computing result does not include a factor that consists both of  $(a_1{+}a_2{+}a_3{+}a_4{+}a_5)$ and $(m_1{+}m_2{+}m_3{+}m_4)$, for example $3(a_1{+}a_2{+}a_3{+}a_4{+}a_5){\rm -}2(m_1{+}m_2{+}m_3{+}m_4)$.
%
\\\\
\newpage \noindent 4. One more theorem\\

Now we present a new theorem.\\\\
\noindent Theorem 2.
\begin{eqnarray}
&& {\rm \ \ \ Let \ } a_1,a_2,\ldots,a_n{\rm \ be\ distinct\ positive\ integers\ such\ that\ }  \nonumber \\
&& 0{<}a_1{<}a_2,\ldots{<}a_n, {\rm \ and }{\ m_1,m_2,\ldots,m_{n-1}}{\rm \ be\ distinct\ positive\ integers}.  \nonumber \\
&& {\rm \ \ \ Now\ if\ any \ two\ distinct\ subsets\ of\ } \{a_1,a_2,\ldots,a_n\}, \nonumber \\
&& \{r_1,r_2,\ldots,r_t\} {\rm \ and \ }\{s_1,s_2,\ldots,s_u\}, {\rm \ hold}  \nonumber \\
&& \ \ \ \ \ \ \ \ \sum_{\rm v=1}^{\rm t} {\rm r}_{\rm v}\neq\sum_{\rm w=1}^{\rm u} {\rm s}_{\rm w} , \\
&& {\rm then\ the\ claim\ of\ the\ grasshopper\ problem\ follows.} \nonumber 
\end{eqnarray}
\begin{proof}
\quad\par
There are totally n! expressions in the form of (3) for degree n.
For any above-mentioned subset if the sum total of each element is equal to one element of $\{m_1,m_2,\ldots,m_{n-1}\}$, then any expression that includes the above-mentioned sum total is equal to 0.

For example if $(a_1{+}a_2{-}m_2)$=0 then

\ \ $f(a_1,a_2,\ldots,a_n)$=${\prod_{l=1}^{\rm n-1} {\prod_{i=1}^{{\rm n}{\rm-1}}}}{((a_1}{+}{a_2}{+}\cdots{+}a_l){-}m_i)$=0

Also there are n!/${}_{\rm n} {\rm C_2}$ expressions, which are in the form of (3), that include $(a_1{+}a_2{-}m_2)$ in.

Whenever the fact that $(x_{\pi(1)}{+}x_{\pi(2)}{+}\cdots{+}x_{\pi(l)}){-}m_i$=0 is found, then the expressions in the form of (3), the values of which are 0, newly increased, in the condition that there is at least one expression that includes $(x_{\pi(1)}{+}x_{\pi(2)}{+}\cdots{+}x_{\pi(l)}){-}m_i$ and its value is not found to be 0 yet. The number of increase is, at most, n!/(${}_{\rm n} {\rm C_l}$). For any l the largest increasing number is n!/(${}_{\rm n} {\rm C_1}$)=(n${\rm -}$1)! , because  ${}_{\rm n} {\rm C_1}$${\leq}$${}_{\rm n} {\rm C_l}$ .

According to the assumption above, the possible largest number of the expressions whose values are 0 is (n${\rm -}$1)!$\times$(n${\rm -}$1).

As a result at least n!${\rm -}$(n${\rm -}$1)!(n${\rm -}$1)$=$(n${\rm -}$1)! expressions remain to be nonzero.

This completes Theorem 2.

\end{proof}
In fact, the condition (21) above is not necessarily guaranteed [4], so we can not apply Theorem 2 easily. 

For the case n=4, there are 10 pairs for which we can not determine one of the two integers are equal to the other, as follows, ($a_1{+}a_2,a_3$),($a_1{+}a_2{+}a_4,a_3{+}a_4$),\\
($a_1{+}a_2,a_4$),($a_1{+}a_2{+}a_3,a_3{+}a_4$),($a_1{+}a_3,a_4$),($a_1{+}a_2{+}a_3,a_2{+}a_4$),\\
($a_2{+}a_3,a_4$),($a_1{+}a_2{+}a_3,a_1{+}a_4$),($a_1{+}a_4,a_2{+}a_3$),($a_1{+}a_2{+}a_3,a_4$).\\

From now on we will prove the claim of the grasshopper problem follows for the case n=4, in spite of the existence of the above-mentioned pairs.\\

Now for the above-mentioned pairs we allocate code names to the event that the two integers are equal to each other, as Table 1.
\begin{table}[h]
\caption{}
\centering
\begin{tabular}{|c|c|}\hline
event & code name \\ \hline
$a_1{+}a_2{=}a_3$ & A1 \\ \hline
$a_1{+}a_2{+}a_4{=}a_3{+}a_4$ & A2 \\ \hline
$a_1{+}a_2{=}a_4$ & B1 \\ \hline
$a_1{+}a_2{+}a_3{=}a_3{+}a_4$ & B2 \\ \hline
$a_1{+}a_3{=}a_4$ & C1 \\ \hline
$a_1{+}a_2{+}a_3{=}a_2{+}a_4$ & C2 \\ \hline
$a_2{+}a_3{=}a_4$ & D1 \\ \hline
$a_1{+}a_2{+}a_3{=}a_1{+}a_4$ & D2 \\ \hline
$a_1{+}a_4{=}a_2{+}a_3$ & E \\ \hline
$a_1{+}a_2{+}a_3{=}a_4$ & F \\ \hline
\end{tabular}
\end{table}
\\\\\\\\\\\\\\\\\\\\\\\\\\

Table 2 shows any pair of 10 events above occurs at the same time. $\bigcirc$ shows the case whenever one of the pair occurs, the other occurs. $\bigtriangleup$ shows the case when there are two partial cases when the other occurs and when the other does not occur, under the condition that one of the pair occurs. $\times$ shows the case whenever one of the pair occurs, the other does not occurs.\\
\begin{table}[h]
\caption{}
\centering
\begin{tabular}{|c|c|c|c|c|c|c|c|c|c|c|}\hline
\ $\backslash$ & A1& A2 & B1 & B2 & C1 & C2 & D1& D2 & E & F \\ \hline
A1 & - &$\bigcirc$&$\times$&$\times$&$\bigtriangleup$&$\bigtriangleup$&$\bigtriangleup$&$\bigtriangleup$&$\bigtriangleup$&$\bigtriangleup$\\ \hline
A2 &$\bigcirc$& - &$\times$&$\times$&$\bigtriangleup$&$\bigtriangleup$&$\bigtriangleup$&$\bigtriangleup$&$\bigtriangleup$&$\bigtriangleup$\\ \hline
B1 & $\times$ & $\times$&-&$\bigcirc$&$\times$&$\times$&$\times$&$\times$&$\bigtriangleup$&$\times$\\ \hline
B2 & $\times$ & $\times$&$\bigcirc$&-&$\times$&$\times$&$\times$&$\times$&$\bigtriangleup$&$\times$\\ \hline
C1 & $\bigtriangleup$& $\bigtriangleup$&$\times$&$\times$&-&$\bigcirc$&$\times$&$\times$&$\bigtriangleup$&$\times$\\ \hline
C2 & $\bigtriangleup$& $\bigtriangleup$&$\times$&$\times$&$\bigcirc$&-&$\times$&$\times$&$\bigtriangleup$&$\times$\\ \hline
D1 & $\bigtriangleup$& $\bigtriangleup$&$\times$&$\times$&$\times$&$\times$&-&$\bigcirc$&$\times$&$\times$\\ \hline
D2 & $\bigtriangleup$& $\bigtriangleup$&$\times$&$\times$&$\times$&$\times$&$\bigcirc$&-&$\times$&$\times$\\ \hline
E & $\bigtriangleup$& $\bigtriangleup$&$\bigtriangleup$&$\bigtriangleup$&$\bigtriangleup$&$\bigtriangleup$&$\times$&$\times$&-&$\times$\\ \hline
F & $\bigtriangleup$& $\bigtriangleup$&$\times$&$\times$&$\times$&$\times$&$\times$&$\times$&$\times$&-\\ \hline
\end{tabular}
\end{table}
\\\\

There are 9 cases how the above-mentioned events occur at a time. For each case we present the combination of occurring events as follows; (A1,A2,C1,C2,E),(A1,A2,D1,D2),(A1,A2,E),(A1,A2,F),(B1,B2,E),(C1,C2,E),\\
(D1,D2),(E),(F).

We assume the contrary, that the above-mentioned 24 polynomials of $a_1$,$a_2$,$a_3$ and $a_4$, made by permutation for the case n=4, are entirely equal to 0.

Table 3 below stores 24 records consisting of a grasshopper's 3 jumping points between the start and the goal which correspond to each expression we have just mentioned above.

\begin{table}[h]
\caption{}
\centering
\begin{tabular}{|c|c|c|c|}\hline
No. & 1st point & 2nd point & 3rd point\\ \hline
1 & $a_1$ & $a_1{+}a_2$ & $a_1{+}a_2{+}a_3$ \\ \hline
2 & $a_1$ & $a_1{+}a_2$ & $a_1{+}a_2{+}a_4$ \\ \hline
3 & $a_1$ & $a_1{+}a_3$ & $a_1{+}a_2{+}a_3$ \\ \hline
4 & $a_1$ & $a_1{+}a_3$ & $a_1{+}a_3{+}a_4$ \\ \hline
5 & $a_1$ & $a_1{+}a_4$ & $a_1{+}a_2{+}a_4$ \\ \hline
6 & $a_1$ & $a_1{+}a_4$ & $a_1{+}a_3{+}a_4$ \\ \hline
7 & $a_2$ & $a_1{+}a_2$ & $a_1{+}a_2{+}a_3$ \\ \hline
8 & $a_2$ & $a_1{+}a_2$ & $a_1{+}a_2{+}a_4$ \\ \hline
9 & $a_2$ & $a_2{+}a_3$ & $a_1{+}a_2{+}a_3$ \\ \hline
10 & $a_2$ & $a_2{+}a_3$ & $a_2{+}a_3{+}a_4$ \\ \hline
11 & $a_2$ & $a_2{+}a_4$ & $a_1{+}a_2{+}a_4$ \\ \hline
12 & $a_2$ & $a_2{+}a_4$ & $a_2{+}a_3{+}a_4$ \\ \hline
13 & $a_3$ & $a_1{+}a_3$ & $a_1{+}a_2{+}a_3$ \\ \hline
14 & $a_3$ & $a_1{+}a_3$ & $a_1{+}a_3{+}a_4$ \\ \hline
15 & $a_3$ & $a_2{+}a_3$ & $a_1{+}a_2{+}a_3$ \\ \hline
16 & $a_3$ & $a_2{+}a_3$ & $a_2{+}a_3{+}a_4$ \\ \hline
17 & $a_3$ & $a_3{+}a_4$ & $a_1{+}a_3{+}a_4$ \\ \hline
18 & $a_3$ & $a_3{+}a_4$ & $a_3{+}a_3{+}a_4$ \\ \hline
19 & $a_4$ & $a_1{+}a_4$ & $a_1{+}a_2{+}a_4$ \\ \hline
20 & $a_4$ & $a_1{+}a_4$ & $a_1{+}a_3{+}a_4$ \\ \hline
21 & $a_4$ & $a_2{+}a_4$ & $a_1{+}a_2{+}a_4$ \\ \hline
22 & $a_4$ & $a_2{+}a_4$ & $a_2{+}a_3{+}a_4$ \\ \hline
23 & $a_4$ & $a_3{+}a_4$ & $a_1{+}a_3{+}a_4$ \\ \hline
24 & $a_4$ & $a_3{+}a_4$ & $a_2{+}a_3{+}a_4$ \\ \hline
\end{tabular}
\end{table}

Now we will arrange Table 3 so that the points of the same value is expressed identically, by using new variables such as $j_1$, $j_2$ and the like.

For the case of  (A1,A2,C1,C2,E).\\

\begin{table}[h]
\caption{}
\centering
\begin{tabular}{|c|c|c|c|}\hline
No. & 1st point & 2nd point & 3rd point\\ \hline
1 & $a_1$ & $j_1$ & $j_4$ \\ \hline
2 & $a_1$ & $j_1$ & $j_2$ \\ \hline
3 & $a_1$ & $j_3$ & $j_4$ \\ \hline
4 & $a_1$ & $j_3$ & $a_1{+}a_3{+}a_4$ \\ \hline
5 & $a_1$ & $j_5$ & $j_2$ \\ \hline
6 & $a_1$ & $j_5$ & $a_1{+}a_3{+}a_4$ \\ \hline
7 & $a_2$ & $j_1$ & $j_4$ \\ \hline
8 & $a_2$ & $j_1$ & $j_2$ \\ \hline
9 & $a_2$ & $j_5$ & $j_4$ \\ \hline
10 & $a_2$ & $j_5$ & $a_2{+}a_3{+}a_4$ \\ \hline
11 & $a_2$ & $j_4$ & $j_2$ \\ \hline
12 & $a_2$ & $j_4$ & $a_2{+}a_3{+}a_4$ \\ \hline
13 & $j_1$ & $j_3$ & $j_4$ \\ \hline
14 & $j_1$ & $j_3$ & $a_1{+}a_3{+}a_4$ \\ \hline
15 & $j_1$ & $j_5$ & $j_4$ \\ \hline
16 & $j_1$ & $j_5$ & $a_2{+}a_3{+}a_4$ \\ \hline
17 & $j_1$ & $j_2$ & $a_1{+}a_3{+}a_4$ \\ \hline
18 & $j_1$ & $j_2$ & $a_2{+}a_3{+}a_4$ \\ \hline
19 & $j_3$ & $j_5$ & $j_2$ \\ \hline
20 & $j_3$ & $j_5$ & $a_1{+}a_3{+}a_4$ \\ \hline
21 & $j_3$ & $j_4$ & $j_2$ \\ \hline
22 & $j_3$ & $j_4$ & $a_2{+}a_3{+}a_4$ \\ \hline
23 & $j_3$ & $j_2$ & $a_1{+}a_3{+}a_4$ \\ \hline
24 & $j_3$ & $j_2$ & $a_2{+}a_3{+}a_4$ \\ \hline
\end{tabular}
\end{table}

\newpage Now we pick up No.6 and No.7, which is convenient for our proof, since these 2 records do not have a point of the same value commonly.

First we assume for 1 of the 2 records its 3 points are equal to $m_1$ or $m_2$ or $m_3$. If for No.6 we assume $a_1$ and $j_5$ and $a_1{+}a_3{+}a_4$ are equal to $m_1$ or $m_2$ or $m_3$, then according to the assumption for No.7 at least one point of $a_2$ and $j_1$ and $j_4$ is equal to $m_1$ or $m_2$ or $m_3$. Then it follows that at least 4 distinct points are equal to $m_1$ or $m_2$ or $m_3$, so the claim can not follow.  If for No.7 we assume $a_2$ and $j_1$ and $j_4$ are equal to $m_1$ or $m_2$ or $m_3$, then according to the assumption for No.6 at least one point of $a_1$ and $j_5$ and $a_1{+}a_3{+}a_4$ is equal to $m_1$ or $m_2$ or $m_3$. Similarly as above, the claim does not follow.

Next we assume for 1 of the 2 records as many as 2 of the 3 points are equal to $m_1$ or $m_2$ or $m_3$. If for No.6 as many as 2 points of $a_1$ and $j_5$ and $a_1{+}a_3{+}a_4$ are equal to $m_1$ or $m_2$ or $m_3$, then for No.7 only 1 point of $a_2$ and $j_1$ and $j_4$ is equal to $m_1$ or $m_2$ or $m_3$ due to the assumption. If the combination of 3 points is ($a_1$,$j_5$,$a_2$), then for No.13 $j_1$ or $j_3$ or $j_4$ is equal to $m_1$ or $m_2$ or $m_3$. It follows that at least 4 distinct points are equal to $m_1$ or $m_2$ or $m_3$, so the claim can not follow. Considering the 9 combinations of 3 points, respectively, No.13 for ($a_1$,$j_5$,$a_2$), No.11 for ($a_1$,$j_5$,$j_1$), No.8 for ($a_1$,$j_5$,$j_4$), No.13 for ($a_1$,$a_1{+}a_3{+}a_4$,$a_2$), No.9 for ($a_1$,$a_1{+}a_3{+}a_4$,$j_1$), No.8 for ($a_1$,$a_1{+}a_3{+}a_4$,$j_4$), No.1 for ($j_5$,$a_1{+}a_3{+}a_4$,$a_2$), No.1 for ($j_5$,$a_1{+}a_3{+}a_4$,$j_1$) and No.2 for ($j_5$,$a_1{+}a_3{+}a_4$,$j_4$) do not make the claim follow. On the contrary, if for No.7 as many as 2 points of $a_2$ and $j_1$ and $j_4$ are equal to $m_1$ or $m_2$ or $m_3$, for No.6 only 1 point of $a_1$ and $j_5$ and $a_1{+}a_3{+}a_4$ is equal to $m_1$ or $m_2$ or $m_3$ due to the assumption. If the combination of 3 points is ($a_1$,$a_2$,$j_1$), then for No.19 $j_3$ or $j_5$ or $j_2$ is equal to $m_1$ or $m_2$ or $m_3$. It follows that at least 4 distinct points are equal to $m_1$ or $m_2$ or $m_3$, so the claim can not follow. Considering the 9 combinations of 3 points, respectively, No.19 for ($a_1$,$a_2$,$j_1$), No.14 for ($a_1$,$a_2$,$j_4$), No.10 for ($a_1$,$j_1$,$j_4$), No.3 for ($j_5$,$a_2$,$j_1$), No.2 for ($j_5$,$a_2$,$j_4$), No.4 for ($j_5$,$j_1$,$j_4$), No.13 for ($a_1{+}a_3{+}a_4$,$a_2$,$j_1$), No.1 for ($a_1{+}a_3{+}a_4$,$a_2$,$j_4$) and No.5 for ($a_1{+}a_3{+}a_4$,$j_1$,$j_4$) do not make the claim follow. 

Last we assume that for both No.6 and No.7 only 1 of the 3 points are equal to $m_1$ or $m_2$ or $m_3$. If the combination of 2 points is ($a_1$,$a_2$), then for No.15 $j_1$ or $j_5$ or $j_4$ is equal to $m_1$ or $m_2$ or $m_3$ and for No.24 $j_3$ or $j_2$ or $a_2{+}a_3{+}a_4$ is equal to $m_1$ or $m_2$ or $m_3$. It follows that at least 4 distinct points are equal to $m_1$ or $m_2$ or $m_3$, so the claim can not follow. Considering the 9 combinations of 2 points, respectively, No.15 and No.24 for ($a_1$,$a_2$), No.10 and No.21 for ($a_1$,$j_1$), No.8 and No.20 for ($a_1$,$j_4$), No.1 and No.21 for ($j_5$,$a_2$), No.4 and No.11 for ($j_5$,$j_1$), No.4 and No.8 for ($j_5$,$j_4$), No.3 and No.18 for ($a_1{+}a_3{+}a_4$,$a_2$), No.3 and No.10 for ($a_1{+}a_3{+}a_4$,$j_1$) and No.2 and No.10 for ($a_1{+}a_3{+}a_4$,$j_4$)  do not make the claim follow.

Therefore for the case of (A1,A2,C1,C2,E), at least 1 of the above-mentioned 24 polynomials of $a_1$,$a_2$,$a_3$ and $a_4$, made by permutation for the case n=4, is not equal to 0, which means that the claim for the grasshopper problem follows. We have succeeded in proving the problem partially.

Similarly for the case of (A1,A2,D1,D2),(A1,A2,E),(A1,A2,F),(B1,B2,E),\\
(C1,C2,E),(D1,D2),(E) and (F), we can also prove that the claim for the grasshopper problem follows.
\\\\
\noindent 5. Proof of the grasshopper problem\\

For perfection we show a proof for the grasshopper problem of ours, we prove it elementarily and inductively.

Here we show the grasshopper problem again.\\

Let $a_1, a_2,\ldots, a_n$ be distinct positive integers and let M be a set of n$-$1 positive integers not containing s=$a_1$+$a_2$$+$$\cdots$$+$$a_n$. A grasshopper is to jump along
the real axis, starting at the point 0 and making n jumps to the right with lengths
$a_1, a_2,\ldots, a_n$ in some order. Prove that the order can be chosen in such a way that
the grasshopper never lands on any point in M.
\begin{proof}
\quad\par
Let A be a set consisting of $a_1, a_2,\ldots, a_n$. Without the loss of generality, we can denote the largest element of A by $a_1$.

And for M=$\{m_1,m_2,\ldots,m_{n-1}\}$ we suppose $m_1{<}m_2{<}\cdots{<}m_{n-1}$.

There are totally 5 cases for the relation of $a_1$ and $m_1$ as follows.

\ \ (a) $a_1{<}m_1$

\ \ (b) $a_1${=}$m_1$

\ \ (c) $a_1{>}m_1$ and $a_1{<}m_{n-1}$ and $a_1{\neq}m_j$

\ \ \ \ \ \ (for any integer j such that 1$\leq$j$\leq$n${\rm -}$1)

\ \ (d) $a_1{>}m_{n-1}$

\ \ (e) $a_1{=}m_j$(for an integer j such that 2$\leq$j$\leq$n${\rm -}$1)

We can prove inductively.

When n=2, A=$\{a_1,a_2\}$ and M=$\{m_1\}$. There are two jumping orders, $(a_1,a_2)$ and $(a_2,a_1)$. According to assumption, $a_1{\neq}m_1$ or $a_2{\neq}m_1$. As a result for at least one of the two orders the claim of this problem follows, that is to say, then a grasshopper can succeed in jumping without blocked.

When n$\leq$k, we assume that the claim of the problem follows  
\\for A=$\{a_1,a_2,\ldots,a_n\}$ and M=$\{m_1,m_2,\ldots,m_{n-1}\}$. In this case, we may regard that any point in M exists between 0 and $a_1$+$a_2$$+$$\cdots$$+$$a_n$, then the claim of this problem still follows.

\noindent ${[}$For the case (a)${]}$

Suppose that when n=k+1 a grasshopper selects a jumping order,\\
$(a_1,a_2,\ldots,a_{k+1})$. Now $a_1{<}m_1$. We consider a series of k jumps ($a_2,\ldots,a_k,\\
a_{k+1}$) starting from $a_1$. For this case, we may regard A=$\{a_2,\ldots,a_k,a_{k+1}\}$ and M=$\{m_1{-}a_1,\ldots,m_k{-}a_1\}$ on the basis of $a_1$. There are ${\rm k}$ points in M between 0 and $a_2$$+$$\cdots$$+$$a_{k+1}$. So the claim of the problem does not follow. Now we temporally omit $m_1{-}a_1$ out of M. Therefore M=$\{m_2{-}a_1,\ldots,m_k{-}a_1\}$. So the claim of the problem follows, in short, there is at least one permutation of $(a_2,\ldots,a_k,a_{k+1})$ that let a grasshopper jump safe without blocked.

If we denote this series of k jumps by ($a_{h2},\ldots,a_{hk},a_{h(k+1)}$), the total series of k+1 jumps ($a_1,a_{h2},\ldots,a_{hk},a_{h(k+1)}$) does not let a grasshopper jump safe, because in fact above-mentioned $m_1$, a point in M, still exists and has a possibility of being landed on by a grasshopper. If not, a grasshopper can jump safe, but if so, there exists an integer l such that 2$\leq$l$\leq$k and 
$a_1{+}a_{h2}{+}\cdots{+}a_{hl}=m_1$. Then by exchanging the first jump for the (l+1)-th jump, we get ($a_{h(l+1)},a_2,\ldots,a_{h(l-1)},a_1,a_{h(l+1)},\ldots,a_{h(k+1)}$) that let a grasshopper jump safe.

As a result the claim of the problem follows.

\noindent ${[}$For the case (b)${]}$

Suppose that when n=k+1 a grasshopper selects a jumping order,\\
$(a_1,a_2,\ldots,a_{k+1})$. Now $a_1{=}m_1$. We consider a series of k jumps ($a_2,\ldots,a_k,\\
a_{k+1}$) starting from $a_1$. For this case, we may regard A=$\{a_2,\ldots,a_k,a_{k+1}\}$ and M=$\{m_2{-}a_1,\ldots,m_k{-}a_1\}$ on the basis of $a_1$. There are ${\rm k}{-}{\rm 1}$ points in M between 0 and $a_2$$+$$\cdots$$+$$a_{k+1}$. So the claim of the problem follows, in short, there is at least one permutation of $(a_2,\ldots,a_k,a_{k+1})$ that let a grasshopper jump safe without blocked.

If we denote this series of k jumps by ($a_{h2},\ldots,a_{hk},a_{h(k+1)}$), the total series of k+1 jumps ($a_1,a_{h2},\ldots,a_{hk},a_{h(k+1)}$) does not let a grasshopper jump safe, because only $a_1$ is a point in M. Then by exchanging the first jump for the second jump, we get ($a_{h2},a_1,\ldots,a_{hk},a_{h(k+1)}$) that let a grasshopper jump safe.

As a result the claim of the problem follows.
 
\noindent ${[}$For the case (c)${]}$

Suppose that when n=k+1 a grasshopper selects a jumping order,\\
$(a_1,a_2,\ldots,a_{k+1})$. Now $m_j{<}a_1{<}m_{j+1}$(for an integer j such that 2$\leq$j$\leq$k). We consider a series of k jumps ($a_2,\ldots,a_k,a_{k+1}$) starting from $a_1$. For this case, we may regard A=$\{a_2,\ldots,a_k,a_{k+1}\}$ and M=$\{m_{j+1}{-}a_1,\ldots,m_k{-}a_1\}$ on the basis of $a_1$. There are ${\rm k}{-}{\rm j}$ points in M between 0 and $a_2$$+$$\cdots$$+$$a_{k+1}$. So the claim of the problem sufficiently follows, in short, there is at least one permutation of $(a_2,\ldots,a_k,a_{k+1})$ that let a grasshopper jump safe without blocked.

Moreover $a_1$ is not any point in M. 

If we denote this series of k jumps by ($a_{h2},\ldots,a_{hk},a_{h(k+1)}$), the total series of k+1 jumps ($a_1,a_{h2},\ldots,a_{hk},a_{h(k+1)}$) let a grasshopper jump safe.

As a result the claim of the problem follows.

\noindent ${[}$For the case (d)${]}$

We easily see the claim of the problem follows.

\noindent ${[}$For the case (e)${]}$

Suppose that when n=k+1 a grasshopper selects a jumping order,\\
$(a_1,a_2,\ldots,a_{k+1})$. Now $a_1{=}m_{j}$(for an integer j such that 2$\leq$j$\leq$k). According to assumption, at least ${\rm k}{-}({\rm j}{-}1)$=${\rm k}{-}{\rm j}{+}1$ elements of a set $\{a_2,\ldots,a_{k+1}\}$ are not equal to any point in M and let $a_g$ be one of its examples.

Now we consider ($a_g,a_1$), which represents the first part sequence of 2 jumps of a sequence of k+1 jumps. The landing point of the first jump is $a_g$, that is not any point in M. And the landing point of the second jump is $a_g{+}a_1$. Note that $m_j{=}a_1{<}a_g{+}a_1{<}a_1{+}a_2{+}\cdots{+}a_{k+1}$ and $m_k{<}a_1{+}a_2{+}\cdots{+}a_{k+1}$. 

There are at most ${\rm k}{-}{\rm j}$ examples that $a_g{+}a_1$ is any point in M. But totally there are at least ${\rm k}{-}{\rm j}{+}1$ examples for $a_g$. Hence a grasshopper succeeds in at least one of the first part sequences of 2 jumps without blocked. Also a grasshopper can jump safe for the second part sequence of k${-}$1 jumps by selecting a suitable jumping order, according to assumption.

As a result the claim of the problem follows. 

\end{proof}
\noindent 6. Discussion and conclusion\\

As we explained in the introduction, it is said that this grasshopper problem can be proved only by elementary and inductive methods(see [1], and we showed above). 

And if they intend to solve by the current method we have shown, there is not perfection yet. 

We can easily assume anti-symmetry of the polynomial function (4). But there is a big drawback, that is to say, 'singularity'. It is not easy to analyze when n is more than 3. 

In short, we are still destined to solve elementarily and deductively, though in most cases, except for 'singularity', a grasshopper succeeds in jumping, judging from (4).

We plan to solve the grasshopper problem by analyzing equations for n's larger than 3 with the aid of Theorem 2.

Last but not least, in the proof of ours above, we do not rely on the condition at all that $a_1, a_2,\ldots, a_n$ and sets of M are integer.  

In short, if the grasshopper problem is as the following,\\

Let $a_1, a_2,\ldots, a_n$ be distinct positive ${\bf numbers}$ and let M be a set of n$-$1 positive ${\bf numbers}$ not containing s=$a_1$+$a_2$$+$$\cdots$$+$$a_n$. A grasshopper is to jump along
the real axis, starting at the point 0 and making n jumps to the right with lengths
$a_1, a_2,\ldots, a_n$ in some order. Prove that the order can be chosen in such a way that
the grasshopper never lands on any point in M. 
\\\\
then the claim of this refined problem still follows.
\\\\\\
\centerline{references}\\
\noindent [1] Geza Kos,6, arxiv, available at https://arxiv.org/abs/1008.2936\\
\noindent [2] N. Alon, Combinatorial Nullstellensatz, Combin. Probab. Comput. 8 (1999) 7–29.\\
\noindent [3] M. Aigner, A Course in Enumeration, Springer-Verlag, Berlin, 2007.\\
\noindent [4] Y. Ieno, vixra, available at http://vixra.org/abs/1910.0283\\
\noindent [5] The Sage Developers. SageMath, the Sage Mathematics Software System (Version 8.9), 2019. https://www.sagemath.org.\\\\

\newpage \noindent Appendix
\\
\\
(20)\
=$(x_1{-}x_2)(x_1{-}x_3)(x_1{-}x_4)(x_1{-}x_5)(x_2{-}x_3)(x_2{-}x_4)(x_2{-}x_5)(x_3{-}x_4)(x_3{-}x_5)(x_4{-}x_5)\\
\times(616x_1^6{+}2440x_1^5x_2{+}4708x_1^4x_2^2{+}5744x_1^3x_2^3{+}4708x_1^2x_2^4{+}2440x_1x_2^5{+}616x_2^6{+}2440x_1^5x_3\\
{+}8332x_1^4x_2x_3{+}13692x_1^3x_2^2x_3{+}13692x_1^2x_2^3x_3{+}8332x_1x_2^4x_3{+}2440x_2^5x_3{+}4708x_1^4x_3^2\\
{+}13692x_1^3x_2x_3^2{+}18436x_1^2x_2^2x_3^2{+}13692x_1x_2^3x_3^2{+}4708x_2^4x_3^2{+}5744x_1^3x_3^3{+}13692x_1^2x_2x_3^3\\
{+}13692x_1x_2^2x_3^3{+}5744x_2^3x_3^3{+}4708x_1^2x_3^4{+}8332x_1x_2x_3^4{+}4708x_2^2x_3^4{+}2440x_1x_3^5{+}2440x_2x_3^5\\
{+}616x_3^6{+}2440x_1^5x_4{+}8332x_1^4x_2x_4{+}13692x_1^3x_2^2x_4{+}13692x_1^2x_2^3x_4{+}8332x_1x_2^4x_4\\
{+}2440x_2^5x_4{+}8332x_1^4x_3x_4{+}24040x_1^3x_2x_3x_4{+}32280x_1^2x_2^2x_3x_4\\
{+}24040x_1x_2^3x_3x_4{+}8332x_2^4x_3x_4{+}13692x_1^3x_3^2x_4{+}32280x_1^2x_2x_3^2x_4{+}32280x_1x_2^2x_3^2x_4\\
{+}13692x_2^3x_3^2x_4{+}13692x_1^2x_3^3x_4{+}24040x_1x_2x_3^3x_4{+}13692x_2^2x_3^3x_4{+}8332x_1x_3^4x_4\\
{+}8332x_2x_3^4x_4{+}2440x_3^5x_4{+}4708x_1^4x_4^2{+}13692x_1^3x_2x_4^2{+}18436x_1^2x_2^2x_4^2{+}13692x_1x_2^3x_4^2\\
{+}4708x_2^4x_4^2{+}13692x_1^3x_3x_4^2{+}32280x_1^2x_2x_3x_4^2{+}32280x_1x_2^2x_3x_4^2{+}13692x_2^3x_3x_4^2\\
{+}18436x_1^2x_3^2x_4^2{+}32280x_1x_2x_3^2x_4^2{+}18436x_2^2x_3^2x_4^2{+}13692x_1x_3^3x_4^2{+}13692x_2x_3^3x_4^2\\
{+}4708x_3^4x_4^2{+}5744x_1^3x_4^3{+}13692x_1^2x_2x_4^3{+}13692x_1x_2^2x_4^3\\
{+}5744x_2^3x_4^3{+}13692x_1^2x_3x_4^3{+}24040x_1x_2x_3x_4^3{+}13692x_2^2x_3x_4^3\\
{+}13692x_1x_3^2x_4^3{+}13692x_2x_3^2x_4^3{+}5744x_3^3x_4^3{+}4708x_1^2x_4^4{+}8332x_1x_2x_4^4\\
{+}4708x_2^2x_4^4{+}8332x_1x_3x_4^4{+}8332x_2x_3x_4^4{+}4708x_3^2x_4^4{+}2440x_1x_4^5{+}2440x_2x_4^5\\
{+}2440x_3x_4^5{+}616x_4^6{+}2440x_1^5x_5{+}8332x_1^4x_2x_5{+}13692x_1^3x_2^2x_5{+}13692x_1^2x_2^3x_5\\
{+}8332x_1x_2^4x_5{+}2440x_2^5x_5{+}8332x_1^4x_3x_5{+}24040x_1^3x_2x_3x_5{+}32280x_1^2x_2^2x_3x_5\\
{+}24040x_1x_2^3x_3x_5{+}8332x_2^4x_3x_5{+}13692x_1^3x_3^2x_5{+}32280x_1^2x_2x_3^2x_5\\
{+}32280x_1x_2^2x_3^2x_5{+}13692x_2^3x_3^2x_5{+}13692x_1^2x_3^3x_5{+}24040x_1x_2x_3^3x_5\\
{+}13692x_2^2x_3^3x_5{+}8332x_1x_3^4x_5{+}8332x_2x_3^4x_5{+}2440x_3^5x_5\\
{+}8332x_1^4x_4x_5{+}24040x_1^3x_2x_4x_5{+}32280x_1^2x_2^2x_4x_5{+}24040x_1x_2^3x_4x_5\\
{+}8332x_2^4x_4x_5{+}24040x_1^3x_3x_4x_5{+}56328x_1^2x_2x_3x_4x_5{+}56328x_1x_2^2x_3x_4x_5\\
{+}24040x_2^3x_3x_4x_5{+}32280x_1^2x_3^2x_4x_5{+}56328x_1x_2x_3^2x_4x_5{+}32280x_2^2x_3^2x_4x_5\\
{+}24040x_1x_3^3x_4x_5{+}24040x_2x_3^3x_4x_5{+}8332x_3^4x_4x_5{+}13692x_1^3x_4^2x_5\\
{+}32280x_1^2x_2x_4^2x_5{+}32280x_1x_2^2x_4^2x_5{+}13692x_2^3x_4^2x_5{+}32280x_1^2x_3x_4^2x_5\\
{+}56328x_1x_2x_3x_4^2x_5{+}32280x_2^2x_3x_4^2x_5{+}32280x_1x_3^2x_4^2x_5{+}32280x_2x_3^2x_4^2x_5\\
{+}13692x_3^3x_4^2x_5{+}13692x_1^2x_4^3x_5{+}24040x_1x_2x_4^3x_5{+}13692x_2^2x_4^3x_5\\
{+}24040x_1x_3x_4^3x_5{+}24040x_2x_3x_4^3x_5{+}13692x_3^2x_4^3x_5{+}8332x_1x_4^4x_5\\
{+}8332x_2x_4^4x_5{+}8332x_3x_4^4x_5{+}2440x_4^5x_5{+}4708x_1^4x_5^2{+}13692x_1^3x_2x_5^2\\
{+}18436x_1^2x_2^2x_5^2{+}13692x_1x_2^3x_5^2{+}4708x_2^4x_5^2{+}13692x_1^3x_3x_5^2\\
{+}32280x_1^2x_2x_3x_5^2{+}32280x_1x_2^2x_3x_5^2{+}13692x_2^3x_3x_5^2{+}18436x_1^2x_3^2x_5^2\\
{+}32280x_1x_2x_3^2x_5^2{+}18436x_2^2x_3^2x_5^2{+}13692x_1x_3^3x_5^2{+}13692x_2x_3^3x_5^2\\
{+}4708x_3^4x_5^2{+}13692x_1^3x_4x_5^2{+}32280x_1^2x_2x_4x_5^2{+}32280x_1x_2^2x_4x_5^2\\
{+}13692x_2^3x_4x_5^2{+}32280x_1^2x_3x_4x_5^2{+}56328x_1x_2x_3x_4x_5^2{+}32280x_2^2x_3x_4x_5^2\\
{+}32280x_1x_3^2x_4x_5^2{+}32280x_2x_3^2x_4x_5^2{+}13692x_3^3x_4x_5^2{+}18436x_1^2x_4^2x_5^2\\
{+}32280x_1x_2x_4^2x_5^2{+}18436x_2^2x_4^2x_5^2{+}32280x_1x_3x_4^2x_5^2{+}32280x_2x_3x_4^2x_5^2\\
{+}18436x_3^2x_4^2x_5^2{+}13692x_1x_4^3x_5^2{+}13692x_2x_4^3x_5^2{+}13692x_3x_4^3x_5^2\\
{+}4708x_4^4x_5^2{+}5744x_1^3x_5^3{+}13692x_1^2x_2x_5^3{+}13692x_1x_2^2x_5^3{+}5744x_2^3x_5^3\\
{+}13692x_1^2x_3x_5^3{+}24040x_1x_2x_3x_5^3{+}13692x_2^2x_3x_5^3{+}13692x_1x_3^2x_5^3\\
{+}13692x_2x_3^2x_5^3{+}5744x_3^3x_5^3{+}13692x_1^2x_4x_5^3{+}24040x_1x_2x_4x_5^3\\
{+}13692x_2^2x_4x_5^3{+}24040x_1x_3x_4x_5^3{+}24040x_2x_3x_4x_5^3{+}13692x_3^2x_4x_5^3\\
{+}13692x_1x_4^2x_5^3{+}13692x_2x_4^2x_5^3{+}13692x_3x_4^2x_5^3{+}5744x_4^3x_5^3\\
{+}4708x_1^2x_5^4{+}8332x_1x_2x_5^4{+}4708x_2^2x_5^4{+}8332x_1x_3x_5^4{+}8332x_2x_3x_5^4\\
{+}4708x_3^2x_5^4{+}8332x_1x_4x_5^4{+}8332x_2x_4x_5^4{+}8332x_3x_4x_5^4{+}4708x_4^2x_5^4\\
{+}2440x_1x_5^5{+}2440x_2x_5^5{+}2440x_3x_5^5{+}2440x_4x_5^5{+}616x_5^6{-}1516x_1^5m_1\\
{-}5404x_1^4x_2m_1{-}9094x_1^3x_2^2m_1{-}9094x_1^2x_2^3m_1{-}5404x_1x_2^4m_1{-}1516x_2^5m_1\\
{-}5404x_1^4x_3m_1{-}16178x_1^3x_2x_3m_1{-}22010x_1^2x_2^2x_3m_1{-}16178x_1x_2^3x_3m_1{-}5404x_2^4x_3m_1\\
{-}9094x_1^3x_3^2m_1{-}22010x_1^2x_2x_3^2m_1{-}22010x_1x_2^2x_3^2m_1{-}9094x_2^3x_3^2m_1\\
{-}9094x_1^2x_3^3m_1{-}16178x_1x_2x_3^3m_1{-}9094x_2^2x_3^3m_1{-}5404x_1x_3^4m_1\\
{-}5404x_2x_3^4m_1{-}1516x_3^5m_1{-}5404x_1^4x_4m_1{-}16178x_1^3x_2x_4m_1{-}22010x_1^2x_2^2x_4m_1{-}16178x_1x_2^3x_4m_1\\
{-}5404x_2^4x_4m_1{-}16178x_1^3x_3x_4m_1{-}38946x_1^2x_2x_3x_4m_1{-}38946x_1x_2^2x_3x_4m_1\\
{-}16178x_2^3x_3x_4m_1{-}22010x_1^2x_3^2x_4m_1{-}38946x_1x_2x_3^2x_4m_1\\
{-}22010x_2^2x_3^2x_4m_1{-}16178x_1x_3^3x_4m_1{-}16178x_2x_3^3x_4m_1\\
{-}5404x_3^4x_4m_1{-}9094x_1^3x_4^2m_1\\
{-}22010x_1^2x_2x_4^2m_1{-}22010x_1x_2^2x_4^2m_1{-}9094x_2^3x_4^2m_1{-}22010x_1^2x_3x_4^2m_1{-}38946x_1x_2x_3x_4^2m_1\\
{-}22010x_2^2x_3x_4^2m_1{-}22010x_1x_3^2x_4^2m_1{-}22010x_2x_3^2x_4^2m_1\\
{-}9094x_3^3x_4^2m_1{-}9094x_1^2x_4^3m_1\\
{-}16178x_1x_2x_4^3m_1{-}9094x_2^2x_4^3m_1{-}16178x_1x_3x_4^3m_1{-}16178x_2x_3x_4^3m_1\\
{-}9094x_3^2x_4^3m_1{-}5404x_1x_4^4m_1\\
{-}5404x_2x_4^4m_1{-}5404x_3x_4^4m_1{-}1516x_4^5m_1{-}5404x_1^4x_5m_1{-}16178x_1^3x_2x_5m_1{-}22010x_1^2x_2^2x_5m_1\\
{-}16178x_1x_2^3x_5m_1{-}5404x_2^4x_5m_1{-}16178x_1^3x_3x_5m_1{-}38946x_1^2x_2x_3x_5m_1\\
{-}38946x_1x_2^2x_3x_5m_1{-}16178x_2^3x_3x_5m_1\\
{-}22010x_1^2x_3^2x_5m_1{-}38946x_1x_2x_3^2x_5m_1{-}22010x_2^2x_3^2x_5m_1\\
{-}16178x_1x_3^3x_5m_1{-}16178x_2x_3^3x_5m_1{-}5404x_3^4x_5m_1\\
{-}16178x_1^3x_4x_5m_1{-}38946x_1^2x_2x_4x_5m_1{-}38946x_1x_2^2x_4x_5m_1{-}16178x_2^3x_4x_5m_1{-}38946x_1^2x_3x_4x_5m_1\\
{-}68700x_1x_2x_3x_4x_5m_1{-}38946x_2^2x_3x_4x_5m_1{-}38946x_1x_3^2x_4x_5m_1\\
{-}38946x_2x_3^2x_4x_5m_1{-}16178x_3^3x_4x_5m_1{-}22010x_1^2x_4^2x_5m_1{-}38946x_1x_2x_4^2x_5m_1{-}22010x_2^2x_4^2x_5m_1\\
{-}38946x_1x_3x_4^2x_5m_1{-}38946x_2x_3x_4^2x_5m_1{-}22010x_3^2x_4^2x_5m_1{-}16178x_1x_4^3x_5m_1{-}16178x_2x_4^3x_5m_1\\
{-}16178x_3x_4^3x_5m_1{-}5404x_4^4x_5m_1{-}9094x_1^3x_5^2m_1{-}22010x_1^2x_2x_5^2m_1{-}22010x_1x_2^2x_5^2m_1\\
{-}9094x_2^3x_5^2m_1{-}22010x_1^2x_3x_5^2m_1{-}38946x_1x_2x_3x_5^2m_1{-}22010x_2^2x_3x_5^2m_1{-}22010x_1x_3^2x_5^2m_1\\
{-}22010x_2x_3^2x_5^2m_1{-}9094x_3^3x_5^2m_1{-}22010x_1^2x_4x_5^2m_1{-}38946x_1x_2x_4x_5^2m_1{-}22010x_2^2x_4x_5^2m_1\\
{-}38946x_1x_3x_4x_5^2m_1{-}38946x_2x_3x_4x_5^2m_1{-}22010x_3^2x_4x_5^2m_1{-}22010x_1x_4^2x_5^2m_1{-}22010x_2x_4^2x_5^2m_1\\
{-}22010x_3x_4^2x_5^2m_1{-}9094x_4^3x_5^2m_1{-}9094x_1^2x_5^3m_1{-}16178x_1x_2x_5^3m_1{-}9094x_2^2x_5^3m_1\\
{-}16178x_1x_3x_5^3m_1{-}16178x_2x_3x_5^3m_1{-}9094x_3^2x_5^3m_1{-}16178x_1x_4x_5^3m_1{-}16178x_2x_4x_5^3m_1\\
{-}16178x_3x_4x_5^3m_1{-}9094x_4^2x_5^3m_1{-}5404x_1x_5^4m_1{-}5404x_2x_5^4m_1{-}5404x_3x_5^4m_1{-}5404x_4x_5^4m_1\\
{-}1516x_5^5m_1{+}1305x_1^4m_1^2{+}4020x_1^3x_2m_1^2{+}5523x_1^2x_2^2m_1^2{+}4020x_1x_2^3m_1^2{+}1305x_2^4m_1^2\\
{+}4020x_1^3x_3m_1^2{+}9883x_1^2x_2x_3m_1^2{+}9883x_1x_2^2x_3m_1^2{+}4020x_2^3x_3m_1^2{+}5523x_1^2x_3^2m_1^2\\
{+}9883x_1x_2x_3^2m_1^2{+}5523x_2^2x_3^2m_1^2{+}4020x_1x_3^3m_1^2{+}4020x_2x_3^3m_1^2{+}1305x_3^4m_1^2\\
{+}4020x_1^3x_4m_1^2{+}9883x_1^2x_2x_4m_1^2{+}9883x_1x_2^2x_4m_1^2{+}4020x_2^3x_4m_1^2{+}9883x_1^2x_3x_4m_1^2\\
{+}17634x_1x_2x_3x_4m_1^2{+}9883x_2^2x_3x_4m_1^2{+}9883x_1x_3^2x_4m_1^2{+}9883x_2x_3^2x_4m_1^2{+}4020x_3^3x_4m_1^2\\
{+}5523x_1^2x_4^2m_1^2{+}9883x_1x_2x_4^2m_1^2{+}5523x_2^2x_4^2m_1^2{+}9883x_1x_3x_4^2m_1^2{+}9883x_2x_3x_4^2m_1^2\\
{+}5523x_3^2x_4^2m_1^2{+}4020x_1x_4^3m_1^2{+}4020x_2x_4^3m_1^2{+}4020x_3x_4^3m_1^2{+}1305x_4^4m_1^2\\
{+}4020x_1^3x_5m_1^2{+}9883x_1^2x_2x_5m_1^2{+}9883x_1x_2^2x_5m_1^2{+}4020x_2^3x_5m_1^2{+}9883x_1^2x_3x_5m_1^2\\
{+}17634x_1x_2x_3x_5m_1^2{+}9883x_2^2x_3x_5m_1^2{+}9883x_1x_3^2x_5m_1^2{+}9883x_2x_3^2x_5m_1^2{+}4020x_3^3x_5m_1^2\\
{+}9883x_1^2x_4x_5m_1^2{+}17634x_1x_2x_4x_5m_1^2{+}9883x_2^2x_4x_5m_1^2{+}17634x_1x_3x_4x_5m_1^2\\
{+}17634x_2x_3x_4x_5m_1^2{+}9883x_3^2x_4x_5m_1^2{+}9883x_1x_4^2x_5m_1^2{+}9883x_2x_4^2x_5m_1^2{+}9883x_3x_4^2x_5m_1^2\\
{+}4020x_4^3x_5m_1^2{+}5523x_1^2x_5^2m_1^2{+}9883x_1x_2x_5^2m_1^2{+}5523x_2^2x_5^2m_1^2{+}9883x_1x_3x_5^2m_1^2\\
{+}9883x_2x_3x_5^2m_1^2{+}5523x_3^2x_5^2m_1^2{+}9883x_1x_4x_5^2m_1^2{+}9883x_2x_4x_5^2m_1^2{+}9883x_3x_4x_5^2m_1^2\\
{+}5523x_4^2x_5^2m_1^2{+}4020x_1x_5^3m_1^2{+}4020x_2x_5^3m_1^2{+}4020x_3x_5^3m_1^2{+}4020x_4x_5^3m_1^2{+}1305x_5^4m_1^2\\
{-}459x_1^3m_1^3{-}1152x_1^2x_2m_1^3{-}1152x_1x_2^2m_1^3{-}459x_2^3m_1^3{-}1152x_1^2x_3m_1^3{-}2079x_1x_2x_3m_1^3\\
{-}1152x_2^2x_3m_1^3{-}1152x_1x_3^2m_1^3{-}1152x_2x_3^2m_1^3{-}459x_3^3m_1^3{-}1152x_1^2x_4m_1^3{-}2079x_1x_2x_4m_1^3\\
{-}1152x_2^2x_4m_1^3{-}2079x_1x_3x_4m_1^3{-}2079x_2x_3x_4m_1^3{-}1152x_3^2x_4m_1^3{-}1152x_1x_4^2m_1^3{-}1152x_2x_4^2m_1^3\\
{-}1152x_3x_4^2m_1^3{-}459x_4^3m_1^3{-}1152x_1^2x_5m_1^3{-}2079x_1x_2x_5m_1^3\\
{-}1152x_2^2x_5m_1^3{-}2079x_1x_3x_5m_1^3{-}2079x_2x_3x_5m_1^3{-}1152x_3^2x_5m_1^3\\
{-}2079x_1x_4x_5m_1^3{-}2079x_2x_4x_5m_1^3{-}2079x_3x_4x_5m_1^3{-}1152x_4^2x_5m_1^3\\
{-}1152x_1x_5^2m_1^3{-}1152x_2x_5^2m_1^3{-}1152x_3x_5^2m_1^3{-}1152x_4x_5^2m_1^3\\
{-}459x_5^3m_1^3{+}54x_1^2m_1^4{+}99x_1x_2m_1^4{+}54x_2^2m_1^4{+}99x_1x_3m_1^4\\
{+}99x_2x_3m_1^4{+}54x_3^2m_1^4{+}99x_1x_4m_1^4{+}99x_2x_4m_1^4\\
{+}99x_3x_4m_1^4{+}54x_4^2m_1^4{+}99x_1x_5m_1^4{+}99x_2x_5m_1^4{+}99x_3x_5m_1^4{+}99x_4x_5m_1^4{+}54x_5^2m_1^4{-}1516x_1^5m_2\\
{-}5404x_1^4x_2m_2{-}9094x_1^3x_2^2m_2{-}9094x_1^2x_2^3m_2{-}5404x_1x_2^4m_2\\
{-}1516x_2^5m_2{-}5404x_1^4x_3m_2{-}16178x_1^3x_2x_3m_2\\
{-}22010x_1^2x_2^2x_3m_2{-}16178x_1x_2^3x_3m_2{-}5404x_2^4x_3m_2{-}9094x_1^3x_3^2m_2\\
{-}22010x_1^2x_2x_3^2m_2{-}22010x_1x_2^2x_3^2m_2\\
{-}9094x_2^3x_3^2m_2{-}9094x_1^2x_3^3m_2{-}16178x_1x_2x_3^3m_2\\
{-}9094x_2^2x_3^3m_2{-}5404x_1x_3^4m_2{-}5404x_2x_3^4m_2{-}1516x_3^5m_2\\
{-}5404x_1^4x_4m_2{-}16178x_1^3x_2x_4m_2{-}22010x_1^2x_2^2x_4m_2\\
{-}16178x_1x_2^3x_4m_2{-}5404x_2^4x_4m_2{-}16178x_1^3x_3x_4m_2\\
{-}38946x_1^2x_2x_3x_4m_2{-}38946x_1x_2^2x_3x_4m_2{-}16178x_2^3x_3x_4m_2\\
{-}22010x_1^2x_3^2x_4m_2{-}38946x_1x_2x_3^2x_4m_2{-}22010x_2^2x_3^2x_4m_2\\
{-}16178x_1x_3^3x_4m_2{-}16178x_2x_3^3x_4m_2{-}5404x_3^4x_4m_2{-}9094x_1^3x_4^2m_2\\
{-}22010x_1^2x_2x_4^2m_2{-}22010x_1x_2^2x_4^2m_2\\
{-}9094x_2^3x_4^2m_2{-}22010x_1^2x_3x_4^2m_2{-}38946x_1x_2x_3x_4^2m_2{-}22010x_2^2x_3x_4^2m_2{-}22010x_1x_3^2x_4^2m_2\\
{-}22010x_2x_3^2x_4^2m_2{-}9094x_3^3x_4^2m_2{-}9094x_1^2x_4^3m_2{-}16178x_1x_2x_4^3m_2\\
{-}9094x_2^2x_4^3m_2{-}16178x_1x_3x_4^3m_2\\
{-}16178x_2x_3x_4^3m_2{-}9094x_3^2x_4^3m_2{-}5404x_1x_4^4m_2\\
{-}5404x_2x_4^4m_2{-}5404x_3x_4^4m_2{-}1516x_4^5m_2{-}5404x_1^4x_5m_2\\
{-}16178x_1^3x_2x_5m_2{-}22010x_1^2x_2^2x_5m_2{-}16178x_1x_2^3x_5m_2{-}5404x_2^4x_5m_2\\
{-}16178x_1^3x_3x_5m_2{-}38946x_1^2x_2x_3x_5m_2\\
{-}38946x_1x_2^2x_3x_5m_2{-}16178x_2^3x_3x_5m_2{-}22010x_1^2x_3^2x_5m_2{-}38946x_1x_2x_3^2x_5m_2{-}22010x_2^2x_3^2x_5m_2\\
{-}16178x_1x_3^3x_5m_2{-}16178x_2x_3^3x_5m_2{-}5404x_3^4x_5m_2{-}16178x_1^3x_4x_5m_2\\
{-}38946x_1^2x_2x_4x_5m_2{-}38946x_1x_2^2x_4x_5m_2\\
{-}16178x_2^3x_4x_5m_2{-}38946x_1^2x_3x_4x_5m_2{-}68700x_1x_2x_3x_4x_5m_2\\
{-}38946x_2^2x_3x_4x_5m_2{-}38946x_1x_3^2x_4x_5m_2\\
{-}38946x_2x_3^2x_4x_5m_2{-}16178x_3^3x_4x_5m_2{-}22010x_1^2x_4^2x_5m_2{-}38946x_1x_2x_4^2x_5m_2{-}22010x_2^2x_4^2x_5m_2\\
{-}38946x_1x_3x_4^2x_5m_2{-}38946x_2x_3x_4^2x_5m_2{-}22010x_3^2x_4^2x_5m_2{-}16178x_1x_4^3x_5m_2\\
{-}16178x_2x_4^3x_5m_2{-}16178x_3x_4^3x_5m_2{-}5404x_4^4x_5m_2{-}9094x_1^3x_5^2m_2{-}22010x_1^2x_2x_5^2m_2\\
{-}22010x_1x_2^2x_5^2m_2{-}9094x_2^3x_5^2m_2{-}22010x_1^2x_3x_5^2m_2{-}38946x_1x_2x_3x_5^2m_2{-}22010x_2^2x_3x_5^2m_2\\
{-}22010x_1x_3^2x_5^2m_2{-}22010x_2x_3^2x_5^2m_2{-}9094x_3^3x_5^2m_2{-}22010x_1^2x_4x_5^2m_2{-}38946x_1x_2x_4x_5^2m_2\\
{-}22010x_2^2x_4x_5^2m_2{-}38946x_1x_3x_4x_5^2m_2{-}38946x_2x_3x_4x_5^2m_2{-}22010x_3^2x_4x_5^2m_2{-}22010x_1x_4^2x_5^2m_2\\
{-}22010x_2x_4^2x_5^2m_2{-}22010x_3x_4^2x_5^2m_2{-}9094x_4^3x_5^2m_2{-}9094x_1^2x_5^3m_2{-}16178x_1x_2x_5^3m_2\\
{-}9094x_2^2x_5^3m_2{-}16178x_1x_3x_5^3m_2{-}16178x_2x_3x_5^3m_2{-}9094x_3^2x_5^3m_2{-}16178x_1x_4x_5^3m_2\\
{-}16178x_2x_4x_5^3m_2{-}16178x_3x_4x_5^3m_2{-}9094x_4^2x_5^3m_2{-}5404x_1x_5^4m_2{-}5404x_2x_5^4m_2\\
{-}5404x_3x_5^4m_2{-}5404x_4x_5^4m_2{-}1516x_5^5m_2{+}3450x_1^4m_1m_2{+}10664x_1^3x_2m_1m_2{+}14674x_1^2x_2^2m_1m_2\\
{+}10664x_1x_2^3m_1m_2{+}3450x_2^4m_1m_2{+}10664x_1^3x_3m_1m_2{+}26298x_1^2x_2x_3m_1m_2{+}26298x_1x_2^2x_3m_1m_2\\
{+}10664x_2^3x_3m_1m_2{+}14674x_1^2x_3^2m_1m_2{+}26298x_1x_2x_3^2m_1m_2{+}14674x_2^2x_3^2m_1m_2{+}10664x_1x_3^3m_1m_2\\
{+}10664x_2x_3^3m_1m_2{+}3450x_3^4m_1m_2{+}10664x_1^3x_4m_1m_2\\
{+}26298x_1^2x_2x_4m_1m_2{+}26298x_1x_2^2x_4m_1m_2\\
{+}10664x_2^3x_4m_1m_2{+}26298x_1^2x_3x_4m_1m_2{+}47004x_1x_2x_3x_4m_1m_2{+}26298x_2^2x_3x_4m_1m_2\\
{+}26298x_1x_3^2x_4m_1m_2{+}26298x_2x_3^2x_4m_1m_2{+}10664x_3^3x_4m_1m_2{+}14674x_1^2x_4^2m_1m_2\\
{+}26298x_1x_2x_4^2m_1m_2{+}14674x_2^2x_4^2m_1m_2{+}26298x_1x_3x_4^2m_1m_2{+}26298x_2x_3x_4^2m_1m_2\\
{+}14674x_3^2x_4^2m_1m_2{+}10664x_1x_4^3m_1m_2{+}10664x_2x_4^3m_1m_2{+}10664x_3x_4^3m_1m_2\\
{+}3450x_4^4m_1m_2{+}10664x_1^3x_5m_1m_2{+}26298x_1^2x_2x_5m_1m_2{+}26298x_1x_2^2x_5m_1m_2{+}10664x_2^3x_5m_1m_2\\
{+}26298x_1^2x_3x_5m_1m_2{+}47004x_1x_2x_3x_5m_1m_2{+}26298x_2^2x_3x_5m_1m_2{+}26298x_1x_3^2x_5m_1m_2\\
{+}26298x_2x_3^2x_5m_1m_2{+}10664x_3^3x_5m_1m_2{+}26298x_1^2x_4x_5m_1m_2{+}47004x_1x_2x_4x_5m_1m_2\\
{+}26298x_2^2x_4x_5m_1m_2{+}47004x_1x_3x_4x_5m_1m_2{+}47004x_2x_3x_4x_5m_1m_2{+}26298x_3^2x_4x_5m_1m_2\\
{+}26298x_1x_4^2x_5m_1m_2{+}26298x_2x_4^2x_5m_1m_2{+}26298x_3x_4^2x_5m_1m_2{+}10664x_4^3x_5m_1m_2\\
{+}14674x_1^2x_5^2m_1m_2{+}26298x_1x_2x_5^2m_1m_2{+}14674x_2^2x_5^2m_1m_2\\
{+}26298x_1x_3x_5^2m_1m_2{+}26298x_2x_3x_5^2m_1m_2{+}14674x_3^2x_5^2m_1m_2\\
{+}26298x_1x_4x_5^2m_1m_2{+}26298x_2x_4x_5^2m_1m_2{+}26298x_3x_4x_5^2m_1m_2{+}14674x_4^2x_5^2m_1m_2\\
{+}10664x_1x_5^3m_1m_2{+}10664x_2x_5^3m_1m_2{+}10664x_3x_5^3m_1m_2{+}10664x_4x_5^3m_1m_2\\
{+}3450x_5^4m_1m_2{-}2663x_1^3m_1^2m_2{-}6694x_1^2x_2m_1^2m_2{-}6694x_1x_2^2m_1^2m_2{-}2663x_2^3m_1^2m_2\\
{-}6694x_1^2x_3m_1^2m_2{-}12093x_1x_2x_3m_1^2m_2{-}6694x_2^2x_3m_1^2m_2{-}6694x_1x_3^2m_1^2m_2\\
{-}6694x_2x_3^2m_1^2m_2{-}2663x_3^3m_1^2m_2{-}6694x_1^2x_4m_1^2m_2{-}12093x_1x_2x_4m_1^2m_2\\
{-}6694x_2^2x_4m_1^2m_2{-}12093x_1x_3x_4m_1^2m_2{-}12093x_2x_3x_4m_1^2m_2{-}6694x_3^2x_4m_1^2m_2\\
{-}6694x_1x_4^2m_1^2m_2{-}6694x_2x_4^2m_1^2m_2{-}6694x_3x_4^2m_1^2m_2\\
{-}2663x_4^3m_1^2m_2{-}6694x_1^2x_5m_1^2m_2{-}12093x_1x_2x_5m_1^2m_2{-}6694x_2^2x_5m_1^2m_2\\
{-}12093x_1x_3x_5m_1^2m_2{-}12093x_2x_3x_5m_1^2m_2{-}6694x_3^2x_5m_1^2m_2\\
{-}12093x_1x_4x_5m_1^2m_2{-}12093x_2x_4x_5m_1^2m_2{-}12093x_3x_4x_5m_1^2m_2{-}6694x_4^2x_5m_1^2m_2\\
{-}6694x_1x_5^2m_1^2m_2{-}6694x_2x_5^2m_1^2m_2{-}6694x_3x_5^2m_1^2m_2\\
{-}6694x_4x_5^2m_1^2m_2{-}2663x_5^3m_1^2m_2{+}798x_1^2m_1^3m_2{+}1458x_1x_2m_1^3m_2{+}798x_2^2m_1^3m_2\\
{+}1458x_1x_3m_1^3m_2{+}1458x_2x_3m_1^3m_2{+}798x_3^2m_1^3m_2{+}1458x_1x_4m_1^3m_2\\
{+}1458x_2x_4m_1^3m_2{+}1458x_3x_4m_1^3m_2{+}798x_4^2m_1^3m_2{+}1458x_1x_5m_1^3m_2\\
{+}1458x_2x_5m_1^3m_2{+}1458x_3x_5m_1^3m_2{+}1458x_4x_5m_1^3m_2{+}798x_5^2m_1^3m_2\\
{-}72x_1m_1^4m_2{-}72x_2m_1^4m_2{-}72x_3m_1^4m_2{-}72x_4m_1^4m_2{-}72x_5m_1^4m_2\\
{+}1305x_1^4m_2^2{+}4020x_1^3x_2m_2^2{+}5523x_1^2x_2^2m_2^2{+}4020x_1x_2^3m_2^2{+}1305x_2^4m_2^2\\
{+}4020x_1^3x_3m_2^2{+}9883x_1^2x_2x_3m_2^2{+}9883x_1x_2^2x_3m_2^2{+}4020x_2^3x_3m_2^2{+}5523x_1^2x_3^2m_2^2\\
{+}9883x_1x_2x_3^2m_2^2{+}5523x_2^2x_3^2m_2^2{+}4020x_1x_3^3m_2^2{+}4020x_2x_3^3m_2^2{+}1305x_3^4m_2^2\\
{+}4020x_1^3x_4m_2^2{+}9883x_1^2x_2x_4m_2^2{+}9883x_1x_2^2x_4m_2^2\\
{+}4020x_2^3x_4m_2^2{+}9883x_1^2x_3x_4m_2^2{+}17634x_1x_2x_3x_4m_2^2{+}9883x_2^2x_3x_4m_2^2\\
{+}9883x_1x_3^2x_4m_2^2{+}9883x_2x_3^2x_4m_2^2{+}4020x_3^3x_4m_2^2{+}5523x_1^2x_4^2m_2^2{+}9883x_1x_2x_4^2m_2^2\\
{+}5523x_2^2x_4^2m_2^2{+}9883x_1x_3x_4^2m_2^2{+}9883x_2x_3x_4^2m_2^2{+}5523x_3^2x_4^2m_2^2{+}4020x_1x_4^3m_2^2\\
{+}4020x_2x_4^3m_2^2{+}4020x_3x_4^3m_2^2{+}1305x_4^4m_2^2{+}4020x_1^3x_5m_2^2{+}9883x_1^2x_2x_5m_2^2\\
{+}9883x_1x_2^2x_5m_2^2{+}4020x_2^3x_5m_2^2{+}9883x_1^2x_3x_5m_2^2{+}17634x_1x_2x_3x_5m_2^2\\
{+}9883x_2^2x_3x_5m_2^2{+}9883x_1x_3^2x_5m_2^2{+}9883x_2x_3^2x_5m_2^2{+}4020x_3^3x_5m_2^2\\
{+}9883x_1^2x_4x_5m_2^2{+}17634x_1x_2x_4x_5m_2^2{+}9883x_2^2x_4x_5m_2^2{+}17634x_1x_3x_4x_5m_2^2\\
{+}17634x_2x_3x_4x_5m_2^2{+}9883x_3^2x_4x_5m_2^2{+}9883x_1x_4^2x_5m_2^2{+}9883x_2x_4^2x_5m_2^2\\
{+}9883x_3x_4^2x_5m_2^2{+}4020x_4^3x_5m_2^2{+}5523x_1^2x_5^2m_2^2{+}9883x_1x_2x_5^2m_2^2{+}5523x_2^2x_5^2m_2^2\\
{+}9883x_1x_3x_5^2m_2^2{+}9883x_2x_3x_5^2m_2^2{+}5523x_3^2x_5^2m_2^2{+}9883x_1x_4x_5^2m_2^2\\
{+}9883x_2x_4x_5^2m_2^2{+}9883x_3x_4x_5^2m_2^2{+}5523x_4^2x_5^2m_2^2{+}4020x_1x_5^3m_2^2\\
{+}4020x_2x_5^3m_2^2{+}4020x_3x_5^3m_2^2{+}4020x_4x_5^3m_2^2{+}1305x_5^4m_2^2{-}2663x_1^3m_1m_2^2\\
{-}6694x_1^2x_2m_1m_2^2{-}6694x_1x_2^2m_1m_2^2{-}2663x_2^3m_1m_2^2{-}6694x_1^2x_3m_1m_2^2\\
{-}12093x_1x_2x_3m_1m_2^2{-}6694x_2^2x_3m_1m_2^2{-}6694x_1x_3^2m_1m_2^2{-}6694x_2x_3^2m_1m_2^2\\
{-}2663x_3^3m_1m_2^2{-}6694x_1^2x_4m_1m_2^2{-}12093x_1x_2x_4m_1m_2^2{-}6694x_2^2x_4m_1m_2^2\\
{-}12093x_1x_3x_4m_1m_2^2{-}12093x_2x_3x_4m_1m_2^2{-}6694x_3^2x_4m_1m_2^2{-}6694x_1x_4^2m_1m_2^2\\
{-}6694x_2x_4^2m_1m_2^2{-}6694x_3x_4^2m_1m_2^2{-}2663x_4^3m_1m_2^2{-}6694x_1^2x_5m_1m_2^2\\
{-}12093x_1x_2x_5m_1m_2^2{-}6694x_2^2x_5m_1m_2^2{-}12093x_1x_3x_5m_1m_2^2{-}12093x_2x_3x_5m_1m_2^2\\
{-}6694x_3^2x_5m_1m_2^2{-}12093x_1x_4x_5m_1m_2^2{-}12093x_2x_4x_5m_1m_2^2{-}12093x_3x_4x_5m_1m_2^2\\
{-}6694x_4^2x_5m_1m_2^2{-}6694x_1x_5^2m_1m_2^2{-}6694x_2x_5^2m_1m_2^2{-}6694x_3x_5^2m_1m_2^2\\
{-}6694x_4x_5^2m_1m_2^2{-}2663x_5^3m_1m_2^2{+}1748x_1^2m_1^2m_2^2{+}3190x_1x_2m_1^2m_2^2\\
{+}1748x_2^2m_1^2m_2^2{+}3190x_1x_3m_1^2m_2^2{+}3190x_2x_3m_1^2m_2^2{+}1748x_3^2m_1^2m_2^2\\
{+}3190x_1x_4m_1^2m_2^2{+}3190x_2x_4m_1^2m_2^2{+}3190x_3x_4m_1^2m_2^2{+}1748x_4^2m_1^2m_2^2\\
{+}3190x_1x_5m_1^2m_2^2{+}3190x_2x_5m_1^2m_2^2{+}3190x_3x_5m_1^2m_2^2{+}3190x_4x_5m_1^2m_2^2\\
{+}1748x_5^2m_1^2m_2^2{-}402x_1m_1^3m_2^2{-}402x_2m_1^3m_2^2{-}402x_3m_1^3m_2^2{-}402x_4m_1^3m_2^2\\
{-}402x_5m_1^3m_2^2{+}21m_1^4m_2^2{-}459x_1^3m_2^3{-}1152x_1^2x_2m_2^3{-}1152x_1x_2^2m_2^3\\
{-}459x_2^3m_2^3{-}1152x_1^2x_3m_2^3{-}2079x_1x_2x_3m_2^3{-}1152x_2^2x_3m_2^3{-}1152x_1x_3^2m_2^3\\
{-}1152x_2x_3^2m_2^3{-}459x_3^3m_2^3{-}1152x_1^2x_4m_2^3{-}2079x_1x_2x_4m_2^3{-}1152x_2^2x_4m_2^3\\
{-}2079x_1x_3x_4m_2^3{-}2079x_2x_3x_4m_2^3{-}1152x_3^2x_4m_2^3{-}1152x_1x_4^2m_2^3{-}1152x_2x_4^2m_2^3\\
{-}1152x_3x_4^2m_2^3{-}459x_4^3m_2^3{-}1152x_1^2x_5m_2^3{-}2079x_1x_2x_5m_2^3{-}1152x_2^2x_5m_2^3\\
{-}2079x_1x_3x_5m_2^3{-}2079x_2x_3x_5m_2^3{-}1152x_3^2x_5m_2^3{-}2079x_1x_4x_5m_2^3{-}2079x_2x_4x_5m_2^3\\
{-}2079x_3x_4x_5m_2^3{-}1152x_4^2x_5m_2^3{-}1152x_1x_5^2m_2^3{-}1152x_2x_5^2m_2^3{-}1152x_3x_5^2m_2^3\\
{-}1152x_4x_5^2m_2^3{-}459x_5^3m_2^3{+}798x_1^2m_1m_2^3{+}1458x_1x_2m_1m_2^3{+}798x_2^2m_1m_2^3\\
{+}1458x_1x_3m_1m_2^3{+}1458x_2x_3m_1m_2^3{+}798x_3^2m_1m_2^3{+}1458x_1x_4m_1m_2^3{+}1458x_2x_4m_1m_2^3\\
{+}1458x_3x_4m_1m_2^3{+}798x_4^2m_1m_2^3{+}1458x_1x_5m_1m_2^3{+}1458x_2x_5m_1m_2^3{+}1458x_3x_5m_1m_2^3\\
{+}1458x_4x_5m_1m_2^3{+}798x_5^2m_1m_2^3{-}402x_1m_1^2m_2^3{-}402x_2m_1^2m_2^3{-}402x_3m_1^2m_2^3\\
{-}402x_4m_1^2m_2^3{-}402x_5m_1^2m_2^3{+}54m_1^3m_2^3{+}54x_1^2m_2^4{+}99x_1x_2m_2^4{+}54x_2^2m_2^4\\
{+}99x_1x_3m_2^4{+}99x_2x_3m_2^4{+}54x_3^2m_2^4{+}99x_1x_4m_2^4{+}99x_2x_4m_2^4{+}99x_3x_4m_2^4\\
{+}54x_4^2m_2^4{+}99x_1x_5m_2^4{+}99x_2x_5m_2^4{+}99x_3x_5m_2^4{+}99x_4x_5m_2^4{+}54x_5^2m_2^4\\
{-}72x_1m_1m_2^4{-}72x_2m_1m_2^4{-}72x_3m_1m_2^4{-}72x_4m_1m_2^4{-}72x_5m_1m_2^4{+}21m_1^2m_2^4\\
{-}1516x_1^5m_3{-}5404x_1^4x_2m_3{-}9094x_1^3x_2^2m_3{-}9094x_1^2x_2^3m_3{-}5404x_1x_2^4m_3\\
{-}1516x_2^5m_3{-}5404x_1^4x_3m_3{-}16178x_1^3x_2x_3m_3\\
{-}22010x_1^2x_2^2x_3m_3{-}16178x_1x_2^3x_3m_3\\
{-}5404x_2^4x_3m_3{-}9094x_1^3x_3^2m_3{-}22010x_1^2x_2x_3^2m_3{-}22010x_1x_2^2x_3^2m_3{-}9094x_2^3x_3^2m_3\\
{-}9094x_1^2x_3^3m_3{-}16178x_1x_2x_3^3m_3{-}9094x_2^2x_3^3m_3{-}5404x_1x_3^4m_3{-}5404x_2x_3^4m_3\\
{-}1516x_3^5m_3{-}5404x_1^4x_4m_3{-}16178x_1^3x_2x_4m_3{-}22010x_1^2x_2^2x_4m_3{-}16178x_1x_2^3x_4m_3\\
{-}5404x_2^4x_4m_3{-}16178x_1^3x_3x_4m_3{-}38946x_1^2x_2x_3x_4m_3{-}38946x_1x_2^2x_3x_4m_3\\
{-}16178x_2^3x_3x_4m_3{-}22010x_1^2x_3^2x_4m_3{-}38946x_1x_2x_3^2x_4m_3{-}22010x_2^2x_3^2x_4m_3\\
{-}16178x_1x_3^3x_4m_3{-}16178x_2x_3^3x_4m_3{-}5404x_3^4x_4m_3{-}9094x_1^3x_4^2m_3\\
{-}22010x_1^2x_2x_4^2m_3{-}22010x_1x_2^2x_4^2m_3\\
{-}9094x_2^3x_4^2m_3{-}22010x_1^2x_3x_4^2m_3{-}38946x_1x_2x_3x_4^2m_3{-}22010x_2^2x_3x_4^2m_3\\
{-}22010x_1x_3^2x_4^2m_3{-}22010x_2x_3^2x_4^2m_3\\
{-}9094x_3^3x_4^2m_3{-}9094x_1^2x_4^3m_3\\
{-}16178x_1x_2x_4^3m_3{-}9094x_2^2x_4^3m_3{-}16178x_1x_3x_4^3m_3{-}16178x_2x_3x_4^3m_3\\
{-}9094x_3^2x_4^3m_3{-}5404x_1x_4^4m_3{-}5404x_2x_4^4m_3{-}5404x_3x_4^4m_3{-}1516x_4^5m_3\\
{-}5404x_1^4x_5m_3{-}16178x_1^3x_2x_5m_3{-}22010x_1^2x_2^2x_5m_3{-}16178x_1x_2^3x_5m_3{-}5404x_2^4x_5m_3\\
{-}16178x_1^3x_3x_5m_3{-}38946x_1^2x_2x_3x_5m_3{-}38946x_1x_2^2x_3x_5m_3{-}16178x_2^3x_3x_5m_3\\
{-}22010x_1^2x_3^2x_5m_3{-}38946x_1x_2x_3^2x_5m_3{-}22010x_2^2x_3^2x_5m_3{-}16178x_1x_3^3x_5m_3\\
{-}16178x_2x_3^3x_5m_3{-}5404x_3^4x_5m_3{-}16178x_1^3x_4x_5m_3{-}38946x_1^2x_2x_4x_5m_3\\
{-}38946x_1x_2^2x_4x_5m_3{-}16178x_2^3x_4x_5m_3{-}38946x_1^2x_3x_4x_5m_3{-}68700x_1x_2x_3x_4x_5m_3\\
{-}38946x_2^2x_3x_4x_5m_3{-}38946x_1x_3^2x_4x_5m_3{-}38946x_2x_3^2x_4x_5m_3{-}16178x_3^3x_4x_5m_3\\
{-}22010x_1^2x_4^2x_5m_3{-}38946x_1x_2x_4^2x_5m_3{-}22010x_2^2x_4^2x_5m_3{-}38946x_1x_3x_4^2x_5m_3\\
{-}38946x_2x_3x_4^2x_5m_3{-}22010x_3^2x_4^2x_5m_3{-}16178x_1x_4^3x_5m_3\\
{-}16178x_2x_4^3x_5m_3{-}16178x_3x_4^3x_5m_3\\
{-}5404x_4^4x_5m_3{-}9094x_1^3x_5^2m_3{-}22010x_1^2x_2x_5^2m_3{-}22010x_1x_2^2x_5^2m_3{-}9094x_2^3x_5^2m_3\\
{-}22010x_1^2x_3x_5^2m_3{-}38946x_1x_2x_3x_5^2m_3{-}22010x_2^2x_3x_5^2m_3{-}22010x_1x_3^2x_5^2m_3\\
{-}22010x_2x_3^2x_5^2m_3{-}9094x_3^3x_5^2m_3{-}22010x_1^2x_4x_5^2m_3{-}38946x_1x_2x_4x_5^2m_3\\
{-}22010x_2^2x_4x_5^2m_3{-}38946x_1x_3x_4x_5^2m_3{-}38946x_2x_3x_4x_5^2m_3{-}22010x_3^2x_4x_5^2m_3\\
{-}22010x_1x_4^2x_5^2m_3{-}22010x_2x_4^2x_5^2m_3{-}22010x_3x_4^2x_5^2m_3{-}9094x_4^3x_5^2m_3\\
{-}9094x_1^2x_5^3m_3{-}16178x_1x_2x_5^3m_3{-}9094x_2^2x_5^3m_3{-}16178x_1x_3x_5^3m_3{-}16178x_2x_3x_5^3m_3\\
{-}9094x_3^2x_5^3m_3{-}16178x_1x_4x_5^3m_3{-}16178x_2x_4x_5^3m_3{-}16178x_3x_4x_5^3m_3{-}9094x_4^2x_5^3m_3\\
{-}5404x_1x_5^4m_3{-}5404x_2x_5^4m_3{-}5404x_3x_5^4m_3{-}5404x_4x_5^4m_3{-}1516x_5^5m_3{+}3450x_1^4m_1m_3\\
{+}10664x_1^3x_2m_1m_3{+}14674x_1^2x_2^2m_1m_3{+}10664x_1x_2^3m_1m_3{+}3450x_2^4m_1m_3{+}10664x_1^3x_3m_1m_3\\
{+}26298x_1^2x_2x_3m_1m_3{+}26298x_1x_2^2x_3m_1m_3{+}10664x_2^3x_3m_1m_3{+}14674x_1^2x_3^2m_1m_3\\
{+}26298x_1x_2x_3^2m_1m_3{+}14674x_2^2x_3^2m_1m_3{+}10664x_1x_3^3m_1m_3{+}10664x_2x_3^3m_1m_3\\
{+}3450x_3^4m_1m_3{+}10664x_1^3x_4m_1m_3{+}26298x_1^2x_2x_4m_1m_3{+}26298x_1x_2^2x_4m_1m_3\\
{+}10664x_2^3x_4m_1m_3{+}26298x_1^2x_3x_4m_1m_3{+}47004x_1x_2x_3x_4m_1m_3{+}26298x_2^2x_3x_4m_1m_3\\
{+}26298x_1x_3^2x_4m_1m_3{+}26298x_2x_3^2x_4m_1m_3{+}10664x_3^3x_4m_1m_3{+}14674x_1^2x_4^2m_1m_3\\
{+}26298x_1x_2x_4^2m_1m_3{+}14674x_2^2x_4^2m_1m_3{+}26298x_1x_3x_4^2m_1m_3{+}26298x_2x_3x_4^2m_1m_3\\
{+}14674x_3^2x_4^2m_1m_3{+}10664x_1x_4^3m_1m_3{+}10664x_2x_4^3m_1m_3{+}10664x_3x_4^3m_1m_3\\
{+}3450x_4^4m_1m_3{+}10664x_1^3x_5m_1m_3{+}26298x_1^2x_2x_5m_1m_3{+}26298x_1x_2^2x_5m_1m_3\\
{+}10664x_2^3x_5m_1m_3{+}26298x_1^2x_3x_5m_1m_3{+}47004x_1x_2x_3x_5m_1m_3{+}26298x_2^2x_3x_5m_1m_3\\
{+}26298x_1x_3^2x_5m_1m_3{+}26298x_2x_3^2x_5m_1m_3{+}10664x_3^3x_5m_1m_3\\
{+}26298x_1^2x_4x_5m_1m_3{+}47004x_1x_2x_4x_5m_1m_3{+}26298x_2^2x_4x_5m_1m_3{+}47004x_1x_3x_4x_5m_1m_3\\
{+}47004x_2x_3x_4x_5m_1m_3{+}26298x_3^2x_4x_5m_1m_3{+}26298x_1x_4^2x_5m_1m_3\\
{+}26298x_2x_4^2x_5m_1m_3{+}26298x_3x_4^2x_5m_1m_3\\
{+}10664x_4^3x_5m_1m_3{+}14674x_1^2x_5^2m_1m_3{+}26298x_1x_2x_5^2m_1m_3\\
{+}14674x_2^2x_5^2m_1m_3{+}26298x_1x_3x_5^2m_1m_3{+}26298x_2x_3x_5^2m_1m_3\\
{+}14674x_3^2x_5^2m_1m_3{+}26298x_1x_4x_5^2m_1m_3{+}26298x_2x_4x_5^2m_1m_3{+}26298x_3x_4x_5^2m_1m_3\\
{+}14674x_4^2x_5^2m_1m_3{+}10664x_1x_5^3m_1m_3{+}10664x_2x_5^3m_1m_3\\
{+}10664x_3x_5^3m_1m_3{+}10664x_4x_5^3m_1m_3\\
{+}3450x_5^4m_1m_3{-}2663x_1^3m_1^2m_3{-}6694x_1^2x_2m_1^2m_3{-}6694x_1x_2^2m_1^2m_3\\
{-}2663x_2^3m_1^2m_3{-}6694x_1^2x_3m_1^2m_3\\
{-}12093x_1x_2x_3m_1^2m_3{-}6694x_2^2x_3m_1^2m_3{-}6694x_1x_3^2m_1^2m_3\\
{-}6694x_2x_3^2m_1^2m_3{-}2663x_3^3m_1^2m_3\\
{-}6694x_1^2x_4m_1^2m_3{-}12093x_1x_2x_4m_1^2m_3{-}6694x_2^2x_4m_1^2m_3{-}\\
12093x_1x_3x_4m_1^2m_3{-}12093x_2x_3x_4m_1^2m_3\\
{-}6694x_3^2x_4m_1^2m_3{-}6694x_1x_4^2m_1^2m_3{-}6694x_2x_4^2m_1^2m_3\\
{-}6694x_3x_4^2m_1^2m_3{-}2663x_4^3m_1^2m_3{-}6694x_1^2x_5m_1^2m_3{-}12093x_1x_2x_5m_1^2m_3\\
{-}6694x_2^2x_5m_1^2m_3{-}12093x_1x_3x_5m_1^2m_3{-}12093x_2x_3x_5m_1^2m_3\\
{-}6694x_3^2x_5m_1^2m_3{-}12093x_1x_4x_5m_1^2m_3{-}12093x_2x_4x_5m_1^2m_3{-}12093x_3x_4x_5m_1^2m_3{-}6694x_4^2x_5m_1^2m_3\\
{-}6694x_1x_5^2m_1^2m_3{-}6694x_2x_5^2m_1^2m_3{-}6694x_3x_5^2m_1^2m_3\\
{-}6694x_4x_5^2m_1^2m_3{-}2663x_5^3m_1^2m_3\\
{+}798x_1^2m_1^3m_3{+}1458x_1x_2m_1^3m_3{+}798x_2^2m_1^3m_3{+}1458x_1x_3m_1^3m_3\\
{+}1458x_2x_3m_1^3m_3{+}798x_3^2m_1^3m_3{+}1458x_1x_4m_1^3m_3{+}1458x_2x_4m_1^3m_3{+}1458x_3x_4m_1^3m_3\\
{+}798x_4^2m_1^3m_3{+}1458x_1x_5m_1^3m_3{+}1458x_2x_5m_1^3m_3{+}1458x_3x_5m_1^3m_3\\
{+}1458x_4x_5m_1^3m_3{+}798x_5^2m_1^3m_3{-}72x_1m_1^4m_3{-}72x_2m_1^4m_3\\
{-}72x_3m_1^4m_3{-}72x_4m_1^4m_3{-}72x_5m_1^4m_3\\
{+}3450x_1^4m_2m_3{+}10664x_1^3x_2m_2m_3{+}14674x_1^2x_2^2m_2m_3{+}10664x_1x_2^3m_2m_3{+}3450x_2^4m_2m_3\\
{+}10664x_1^3x_3m_2m_3{+}26298x_1^2x_2x_3m_2m_3{+}26298x_1x_2^2x_3m_2m_3\\
{+}10664x_2^3x_3m_2m_3{+}14674x_1^2x_3^2m_2m_3{+}26298x_1x_2x_3^2m_2m_3{+}14674x_2^2x_3^2m_2m_3\\
{+}10664x_1x_3^3m_2m_3{+}10664x_2x_3^3m_2m_3{+}3450x_3^4m_2m_3{+}10664x_1^3x_4m_2m_3\\
{+}26298x_1^2x_2x_4m_2m_3{+}26298x_1x_2^2x_4m_2m_3{+}10664x_2^3x_4m_2m_3{+}26298x_1^2x_3x_4m_2m_3\\
{+}47004x_1x_2x_3x_4m_2m_3{+}26298x_2^2x_3x_4m_2m_3{+}26298x_1x_3^2x_4m_2m_3\\
{+}26298x_2x_3^2x_4m_2m_3{+}10664x_3^3x_4m_2m_3\\
{+}14674x_1^2x_4^2m_2m_3{+}26298x_1x_2x_4^2m_2m_3{+}14674x_2^2x_4^2m_2m_3\\
{+}26298x_1x_3x_4^2m_2m_3{+}26298x_2x_3x_4^2m_2m_3\\
{+}14674x_3^2x_4^2m_2m_3{+}10664x_1x_4^3m_2m_3{+}10664x_2x_4^3m_2m_3\\
{+}10664x_3x_4^3m_2m_3{+}3450x_4^4m_2m_3\\
{+}10664x_1^3x_5m_2m_3{+}26298x_1^2x_2x_5m_2m_3{+}26298x_1x_2^2x_5m_2m_3\\
{+}10664x_2^3x_5m_2m_3{+}26298x_1^2x_3x_5m_2m_3\\
{+}47004x_1x_2x_3x_5m_2m_3{+}26298x_2^2x_3x_5m_2m_3{+}26298x_1x_3^2x_5m_2m_3\\
{+}26298x_2x_3^2x_5m_2m_3{+}10664x_3^3x_5m_2m_3\\
{+}26298x_1^2x_4x_5m_2m_3{+}47004x_1x_2x_4x_5m_2m_3{+}26298x_2^2x_4x_5m_2m_3{+}47004x_1x_3x_4x_5m_2m_3\\
{+}47004x_2x_3x_4x_5m_2m_3{+}26298x_3^2x_4x_5m_2m_3\\
{+}26298x_1x_4^2x_5m_2m_3{+}26298x_2x_4^2x_5m_2m_3\\
{+}26298x_3x_4^2x_5m_2m_3{+}10664x_4^3x_5m_2m_3{+}14674x_1^2x_5^2m_2m_3{+}26298x_1x_2x_5^2m_2m_3\\
{+}14674x_2^2x_5^2m_2m_3{+}26298x_1x_3x_5^2m_2m_3{+}26298x_2x_3x_5^2m_2m_3{+}14674x_3^2x_5^2m_2m_3\\
{+}26298x_1x_4x_5^2m_2m_3{+}26298x_2x_4x_5^2m_2m_3{+}26298x_3x_4x_5^2m_2m_3{+}14674x_4^2x_5^2m_2m_3\\
{+}10664x_1x_5^3m_2m_3{+}10664x_2x_5^3m_2m_3{+}10664x_3x_5^3m_2m_3{+}10664x_4x_5^3m_2m_3\\
{+}3450x_5^4m_2m_3{-}6912x_1^3m_1m_2m_3{-}17390x_1^2x_2m_1m_2m_3{-}17390x_1x_2^2m_1m_2m_3\\
{-}6912x_2^3m_1m_2m_3{-}17390x_1^2x_3m_1m_2m_3{-}31434x_1x_2x_3m_1m_2m_3{-}17390x_2^2x_3m_1m_2m_3\\
{-}17390x_1x_3^2m_1m_2m_3{-}17390x_2x_3^2m_1m_2m_3{-}6912x_3^3m_1m_2m_3\\
{-}17390x_1^2x_4m_1m_2m_3{-}31434x_1x_2x_4m_1m_2m_3{-}17390x_2^2x_4m_1m_2m_3{-}31434x_1x_3x_4m_1m_2m_3\\
{-}31434x_2x_3x_4m_1m_2m_3{-}17390x_3^2x_4m_1m_2m_3{-}17390x_1x_4^2m_1m_2m_3{-}17390x_2x_4^2m_1m_2m_3\\
{-}17390x_3x_4^2m_1m_2m_3{-}6912x_4^3m_1m_2m_3{-}17390x_1^2x_5m_1m_2m_3{-}31434x_1x_2x_5m_1m_2m_3\\
{-}17390x_2^2x_5m_1m_2m_3{-}31434x_1x_3x_5m_1m_2m_3{-}31434x_2x_3x_5m_1m_2m_3\\
{-}17390x_3^2x_5m_1m_2m_3{-}31434x_1x_4x_5m_1m_2m_3{-}31434x_2x_4x_5m_1m_2m_3\\
{-}31434x_3x_4x_5m_1m_2m_3{-}17390x_4^2x_5m_1m_2m_3{-}17390x_1x_5^2m_1m_2m_3\\
{-}17390x_2x_5^2m_1m_2m_3{-}17390x_3x_5^2m_1m_2m_3\\
{-}17390x_4x_5^2m_1m_2m_3{-}6912x_5^3m_1m_2m_3{+}4452x_1^2m_1^2m_2m_3{+}8122x_1x_2m_1^2m_2m_3\\
{+}4452x_2^2m_1^2m_2m_3{+}8122x_1x_3m_1^2m_2m_3{+}8122x_2x_3m_1^2m_2m_3{+}4452x_3^2m_1^2m_2m_3\\
{+}8122x_1x_4m_1^2m_2m_3{+}8122x_2x_4m_1^2m_2m_3{+}8122x_3x_4m_1^2m_2m_3{+}4452x_4^2m_1^2m_2m_3\\
{+}8122x_1x_5m_1^2m_2m_3{+}8122x_2x_5m_1^2m_2m_3{+}8122x_3x_5m_1^2m_2m_3{+}8122x_4x_5m_1^2m_2m_3\\
{+}4452x_5^2m_1^2m_2m_3{-}1002x_1m_1^3m_2m_3{-}1002x_2m_1^3m_2m_3{-}1002x_3m_1^3m_2m_3\\
{-}1002x_4m_1^3m_2m_3{-}1002x_5m_1^3m_2m_3{+}51m_1^4m_2m_3{-}2663x_1^3m_2^2m_3{-}6694x_1^2x_2m_2^2m_3\\
{-}6694x_1x_2^2m_2^2m_3{-}2663x_2^3m_2^2m_3{-}6694x_1^2x_3m_2^2m_3{-}12093x_1x_2x_3m_2^2m_3\\
{-}6694x_2^2x_3m_2^2m_3{-}6694x_1x_3^2m_2^2m_3{-}6694x_2x_3^2m_2^2m_3{-}2663x_3^3m_2^2m_3\\
{-}6694x_1^2x_4m_2^2m_3{-}12093x_1x_2x_4m_2^2m_3{-}6694x_2^2x_4m_2^2m_3{-}12093x_1x_3x_4m_2^2m_3\\
{-}12093x_2x_3x_4m_2^2m_3{-}6694x_3^2x_4m_2^2m_3{-}6694x_1x_4^2m_2^2m_3{-}6694x_2x_4^2m_2^2m_3\\
{-}6694x_3x_4^2m_2^2m_3{-}2663x_4^3m_2^2m_3{-}6694x_1^2x_5m_2^2m_3{-}12093x_1x_2x_5m_2^2m_3\\
{-}6694x_2^2x_5m_2^2m_3{-}12093x_1x_3x_5m_2^2m_3{-}12093x_2x_3x_5m_2^2m_3\\
{-}6694x_3^2x_5m_2^2m_3{-}12093x_1x_4x_5m_2^2m_3{-}12093x_2x_4x_5m_2^2m_3{-}12093x_3x_4x_5m_2^2m_3\\
{-}6694x_4^2x_5m_2^2m_3{-}6694x_1x_5^2m_2^2m_3{-}6694x_2x_5^2m_2^2m_3{-}6694x_3x_5^2m_2^2m_3\\
{-}6694x_4x_5^2m_2^2m_3{-}2663x_5^3m_2^2m_3{+}4452x_1^2m_1m_2^2m_3{+}8122x_1x_2m_1m_2^2m_3\\
{+}4452x_2^2m_1m_2^2m_3{+}8122x_1x_3m_1m_2^2m_3{+}8122x_2x_3m_1m_2^2m_3\\
{+}4452x_3^2m_1m_2^2m_3{+}8122x_1x_4m_1m_2^2m_3{+}8122x_2x_4m_1m_2^2m_3{+}8122x_3x_4m_1m_2^2m_3\\
{+}4452x_4^2m_1m_2^2m_3{+}8122x_1x_5m_1m_2^2m_3{+}8122x_2x_5m_1m_2^2m_3{+}8122x_3x_5m_1m_2^2m_3\\
{+}8122x_4x_5m_1m_2^2m_3{+}4452x_5^2m_1m_2^2m_3{-}2158x_1m_1^2m_2^2m_3\\
{-}2158x_2m_1^2m_2^2m_3{-}2158x_3m_1^2m_2^2m_3{-}2158x_4m_1^2m_2^2m_3\\
{-}2158x_5m_1^2m_2^2m_3{+}279m_1^3m_2^2m_3{+}798x_1^2m_2^3m_3{+}1458x_1x_2m_2^3m_3{+}798x_2^2m_2^3m_3\\
{+}1458x_1x_3m_2^3m_3{+}1458x_2x_3m_2^3m_3{+}798x_3^2m_2^3m_3{+}1458x_1x_4m_2^3m_3{+}1458x_2x_4m_2^3m_3\\
{+}1458x_3x_4m_2^3m_3{+}798x_4^2m_2^3m_3{+}1458x_1x_5m_2^3m_3{+}1458x_2x_5m_2^3m_3\\
{+}1458x_3x_5m_2^3m_3{+}1458x_4x_5m_2^3m_3\\
{+}798x_5^2m_2^3m_3{-}1002x_1m_1m_2^3m_3{-}1002x_2m_1m_2^3m_3{-}1002x_3m_1m_2^3m_3\\
{-}1002x_4m_1m_2^3m_3{-}1002x_5m_1m_2^3m_3\\
{+}279m_1^2m_2^3m_3{-}72x_1m_2^4m_3{-}72x_2m_2^4m_3{-}72x_3m_2^4m_3{-}72x_4m_2^4m_3{-}72x_5m_2^4m_3\\
{+}51m_1m_2^4m_3{+}1305x_1^4m_3^2{+}4020x_1^3x_2m_3^2{+}5523x_1^2x_2^2m_3^2{+}4020x_1x_2^3m_3^2\\
{+}1305x_2^4m_3^2{+}4020x_1^3x_3m_3^2{+}9883x_1^2x_2x_3m_3^2{+}9883x_1x_2^2x_3m_3^2{+}4020x_2^3x_3m_3^2\\
{+}5523x_1^2x_3^2m_3^2{+}9883x_1x_2x_3^2m_3^2{+}5523x_2^2x_3^2m_3^2{+}4020x_1x_3^3m_3^2{+}4020x_2x_3^3m_3^2\\
{+}1305x_3^4m_3^2{+}4020x_1^3x_4m_3^2\\
{+}9883x_1^2x_2x_4m_3^2{+}9883x_1x_2^2x_4m_3^2{+}4020x_2^3x_4m_3^2{+}9883x_1^2x_3x_4m_3^2\\
{+}17634x_1x_2x_3x_4m_3^2{+}9883x_2^2x_3x_4m_3^2{+}9883x_1x_3^2x_4m_3^2{+}9883x_2x_3^2x_4m_3^2\\
{+}4020x_3^3x_4m_3^2{+}5523x_1^2x_4^2m_3^2{+}9883x_1x_2x_4^2m_3^2\\
{+}5523x_2^2x_4^2m_3^2{+}9883x_1x_3x_4^2m_3^2{+}9883x_2x_3x_4^2m_3^2{+}5523x_3^2x_4^2m_3^2\\
{+}4020x_1x_4^3m_3^2{+}4020x_2x_4^3m_3^2{+}4020x_3x_4^3m_3^2{+}1305x_4^4m_3^2{+}4020x_1^3x_5m_3^2\\
{+}9883x_1^2x_2x_5m_3^2{+}9883x_1x_2^2x_5m_3^2\\
{+}4020x_2^3x_5m_3^2{+}9883x_1^2x_3x_5m_3^2{+}17634x_1x_2x_3x_5m_3^2{+}9883x_2^2x_3x_5m_3^2\\
{+}9883x_1x_3^2x_5m_3^2{+}9883x_2x_3^2x_5m_3^2{+}4020x_3^3x_5m_3^2{+}9883x_1^2x_4x_5m_3^2\\
{+}17634x_1x_2x_4x_5m_3^2{+}9883x_2^2x_4x_5m_3^2{+}17634x_1x_3x_4x_5m_3^2{+}17634x_2x_3x_4x_5m_3^2\\
{+}9883x_3^2x_4x_5m_3^2{+}9883x_1x_4^2x_5m_3^2{+}9883x_2x_4^2x_5m_3^2{+}9883x_3x_4^2x_5m_3^2\\
{+}4020x_4^3x_5m_3^2{+}5523x_1^2x_5^2m_3^2{+}9883x_1x_2x_5^2m_3^2{+}5523x_2^2x_5^2m_3^2\\
{+}9883x_1x_3x_5^2m_3^2{+}9883x_2x_3x_5^2m_3^2{+}5523x_3^2x_5^2m_3^2{+}9883x_1x_4x_5^2m_3^2\\
{+}9883x_2x_4x_5^2m_3^2{+}9883x_3x_4x_5^2m_3^2{+}5523x_4^2x_5^2m_3^2\\
{+}4020x_1x_5^3m_3^2{+}4020x_2x_5^3m_3^2{+}4020x_3x_5^3m_3^2{+}4020x_4x_5^3m_3^2{+}1305x_5^4m_3^2\\
{-}2663x_1^3m_1m_3^2{-}6694x_1^2x_2m_1m_3^2{-}6694x_1x_2^2m_1m_3^2{-}2663x_2^3m_1m_3^2\\
{-}6694x_1^2x_3m_1m_3^2{-}12093x_1x_2x_3m_1m_3^2\\
{-}6694x_2^2x_3m_1m_3^2{-}6694x_1x_3^2m_1m_3^2{-}6694x_2x_3^2m_1m_3^2{-}2663x_3^3m_1m_3^2\\
{-}6694x_1^2x_4m_1m_3^2{-}12093x_1x_2x_4m_1m_3^2{-}6694x_2^2x_4m_1m_3^2{-}12093x_1x_3x_4m_1m_3^2\\
{-}12093x_2x_3x_4m_1m_3^2{-}6694x_3^2x_4m_1m_3^2{-}6694x_1x_4^2m_1m_3^2{-}6694x_2x_4^2m_1m_3^2\\
{-}6694x_3x_4^2m_1m_3^2{-}2663x_4^3m_1m_3^2{-}6694x_1^2x_5m_1m_3^2\\
{-}12093x_1x_2x_5m_1m_3^2{-}6694x_2^2x_5m_1m_3^2{-}12093x_1x_3x_5m_1m_3^2{-}12093x_2x_3x_5m_1m_3^2\\
{-}6694x_3^2x_5m_1m_3^2{-}12093x_1x_4x_5m_1m_3^2{-}12093x_2x_4x_5m_1m_3^2{-}12093x_3x_4x_5m_1m_3^2\\
{-}6694x_4^2x_5m_1m_3^2{-}6694x_1x_5^2m_1m_3^2{-}6694x_2x_5^2m_1m_3^2\\
{-}6694x_3x_5^2m_1m_3^2{-}6694x_4x_5^2m_1m_3^2\\
{-}2663x_5^3m_1m_3^2{+}1748x_1^2m_1^2m_3^2{+}3190x_1x_2m_1^2m_3^2{+}1748x_2^2m_1^2m_3^2\\
{+}3190x_1x_3m_1^2m_3^2{+}3190x_2x_3m_1^2m_3^2{+}1748x_3^2m_1^2m_3^2{+}3190x_1x_4m_1^2m_3^2\\
{+}3190x_2x_4m_1^2m_3^2{+}3190x_3x_4m_1^2m_3^2{+}1748x_4^2m_1^2m_3^2{+}3190x_1x_5m_1^2m_3^2\\
{+}3190x_2x_5m_1^2m_3^2{+}3190x_3x_5m_1^2m_3^2{+}3190x_4x_5m_1^2m_3^2{+}1748x_5^2m_1^2m_3^2\\
{-}402x_1m_1^3m_3^2{-}402x_2m_1^3m_3^2{-}402x_3m_1^3m_3^2{-}402x_4m_1^3m_3^2{-}402x_5m_1^3m_3^2\\
{+}21m_1^4m_3^2{-}2663x_1^3m_2m_3^2{-}6694x_1^2x_2m_2m_3^2{-}6694x_1x_2^2m_2m_3^2{-}2663x_2^3m_2m_3^2\\
{-}6694x_1^2x_3m_2m_3^2{-}12093x_1x_2x_3m_2m_3^2{-}6694x_2^2x_3m_2m_3^2{-}6694x_1x_3^2m_2m_3^2\\
{-}6694x_2x_3^2m_2m_3^2{-}2663x_3^3m_2m_3^2{-}6694x_1^2x_4m_2m_3^2{-}12093x_1x_2x_4m_2m_3^2\\
{-}6694x_2^2x_4m_2m_3^2{-}12093x_1x_3x_4m_2m_3^2{-}12093x_2x_3x_4m_2m_3^2{-}6694x_3^2x_4m_2m_3^2\\
{-}6694x_1x_4^2m_2m_3^2{-}6694x_2x_4^2m_2m_3^2{-}6694x_3x_4^2m_2m_3^2{-}2663x_4^3m_2m_3^2\\
{-}6694x_1^2x_5m_2m_3^2{-}12093x_1x_2x_5m_2m_3^2{-}6694x_2^2x_5m_2m_3^2{-}12093x_1x_3x_5m_2m_3^2\\
{-}12093x_2x_3x_5m_2m_3^2{-}6694x_3^2x_5m_2m_3^2{-}12093x_1x_4x_5m_2m_3^2{-}12093x_2x_4x_5m_2m_3^2\\
{-}12093x_3x_4x_5m_2m_3^2{-}6694x_4^2x_5m_2m_3^2{-}6694x_1x_5^2m_2m_3^2{-}6694x_2x_5^2m_2m_3^2\\
{-}6694x_3x_5^2m_2m_3^2{-}6694x_4x_5^2m_2m_3^2{-}2663x_5^3m_2m_3^2{+}4452x_1^2m_1m_2m_3^2\\
{+}8122x_1x_2m_1m_2m_3^2{+}4452x_2^2m_1m_2m_3^2{+}8122x_1x_3m_1m_2m_3^2{+}8122x_2x_3m_1m_2m_3^2\\
{+}4452x_3^2m_1m_2m_3^2{+}8122x_1x_4m_1m_2m_3^2{+}8122x_2x_4m_1m_2m_3^2{+}8122x_3x_4m_1m_2m_3^2\\
{+}4452x_4^2m_1m_2m_3^2{+}8122x_1x_5m_1m_2m_3^2{+}8122x_2x_5m_1m_2m_3^2{+}8122x_3x_5m_1m_2m_3^2\\
{+}8122x_4x_5m_1m_2m_3^2{+}4452x_5^2m_1m_2m_3^2{-}2158x_1m_1^2m_2m_3^2{-}2158x_2m_1^2m_2m_3^2\\
{-}2158x_3m_1^2m_2m_3^2{-}2158x_4m_1^2m_2m_3^2{-}2158x_5m_1^2m_2m_3^2{+}279m_1^3m_2m_3^2\\
{+}1748x_1^2m_2^2m_3^2{+}3190x_1x_2m_2^2m_3^2{+}1748x_2^2m_2^2m_3^2{+}3190x_1x_3m_2^2m_3^2\\
{+}3190x_2x_3m_2^2m_3^2{+}1748x_3^2m_2^2m_3^2{+}3190x_1x_4m_2^2m_3^2{+}3190x_2x_4m_2^2m_3^2\\
{+}3190x_3x_4m_2^2m_3^2{+}1748x_4^2m_2^2m_3^2{+}3190x_1x_5m_2^2m_3^2{+}3190x_2x_5m_2^2m_3^2\\
{+}3190x_3x_5m_2^2m_3^2{+}3190x_4x_5m_2^2m_3^2{+}1748x_5^2m_2^2m_3^2{-}2158x_1m_1m_2^2m_3^2\\
{-}2158x_2m_1m_2^2m_3^2{-}2158x_3m_1m_2^2m_3^2{-}2158x_4m_1m_2^2m_3^2{-}2158x_5m_1m_2^2m_3^2\\
{+}593m_1^2m_2^2m_3^2{-}402x_1m_2^3m_3^2{-}402x_2m_2^3m_3^2{-}402x_3m_2^3m_3^2{-}402x_4m_2^3m_3^2\\
{-}402x_5m_2^3m_3^2{+}279m_1m_2^3m_3^2{+}21m_2^4m_3^2{-}459x_1^3m_3^3{-}1152x_1^2x_2m_3^3\\
{-}1152x_1x_2^2m_3^3{-}459x_2^3m_3^3{-}1152x_1^2x_3m_3^3{-}2079x_1x_2x_3m_3^3{-}1152x_2^2x_3m_3^3\\
{-}1152x_1x_3^2m_3^3{-}1152x_2x_3^2m_3^3{-}459x_3^3m_3^3{-}1152x_1^2x_4m_3^3{-}2079x_1x_2x_4m_3^3\\
{-}1152x_2^2x_4m_3^3{-}2079x_1x_3x_4m_3^3{-}2079x_2x_3x_4m_3^3\\
{-}1152x_3^2x_4m_3^3{-}1152x_1x_4^2m_3^3{-}1152x_2x_4^2m_3^3\\
{-}1152x_3x_4^2m_3^3{-}459x_4^3m_3^3{-}1152x_1^2x_5m_3^3{-}2079x_1x_2x_5m_3^3{-}1152x_2^2x_5m_3^3\\
{-}2079x_1x_3x_5m_3^3{-}2079x_2x_3x_5m_3^3{-}1152x_3^2x_5m_3^3{-}2079x_1x_4x_5m_3^3{-}2079x_2x_4x_5m_3^3\\
{-}2079x_3x_4x_5m_3^3{-}1152x_4^2x_5m_3^3{-}1152x_1x_5^2m_3^3{-}1152x_2x_5^2m_3^3{-}1152x_3x_5^2m_3^3\\
{-}1152x_4x_5^2m_3^3{-}459x_5^3m_3^3\\
{+}798x_1^2m_1m_3^3{+}1458x_1x_2m_1m_3^3{+}798x_2^2m_1m_3^3{+}1458x_1x_3m_1m_3^3{+}1458x_2x_3m_1m_3^3\\
{+}798x_3^2m_1m_3^3{+}1458x_1x_4m_1m_3^3{+}1458x_2x_4m_1m_3^3{+}1458x_3x_4m_1m_3^3\\
{+}798x_4^2m_1m_3^3{+}1458x_1x_5m_1m_3^3{+}1458x_2x_5m_1m_3^3{+}1458x_3x_5m_1m_3^3\\
{+}1458x_4x_5m_1m_3^3{+}798x_5^2m_1m_3^3\\
{-}402x_1m_1^2m_3^3{-}402x_2m_1^2m_3^3{-}402x_3m_1^2m_3^3{-}402x_4m_1^2m_3^3{-}402x_5m_1^2m_3^3\\
{+}54m_1^3m_3^3{+}798x_1^2m_2m_3^3{+}1458x_1x_2m_2m_3^3{+}798x_2^2m_2m_3^3{+}1458x_1x_3m_2m_3^3\\
{+}1458x_2x_3m_2m_3^3{+}798x_3^2m_2m_3^3{+}1458x_1x_4m_2m_3^3{+}1458x_2x_4m_2m_3^3{+}1458x_3x_4m_2m_3^3\\
{+}798x_4^2m_2m_3^3{+}1458x_1x_5m_2m_3^3{+}1458x_2x_5m_2m_3^3\\
{+}1458x_3x_5m_2m_3^3{+}1458x_4x_5m_2m_3^3{+}798x_5^2m_2m_3^3\\
{-}1002x_1m_1m_2m_3^3{-}1002x_2m_1m_2m_3^3{-}1002x_3m_1m_2m_3^3{-}1002x_4m_1m_2m_3^3\\
{-}1002x_5m_1m_2m_3^3{+}279m_1^2m_2m_3^3{-}402x_1m_2^2m_3^3{-}402x_2m_2^2m_3^3{-}402x_3m_2^2m_3^3\\
{-}402x_4m_2^2m_3^3{-}402x_5m_2^2m_3^3{+}279m_1m_2^2m_3^3{+}54m_2^3m_3^3{+}54x_1^2m_3^4{+}99x_1x_2m_3^4\\
{+}54x_2^2m_3^4{+}99x_1x_3m_3^4{+}99x_2x_3m_3^4{+}54x_3^2m_3^4{+}99x_1x_4m_3^4\\
{+}99x_2x_4m_3^4{+}99x_3x_4m_3^4{+}54x_4^2m_3^4\\
{+}99x_1x_5m_3^4{+}99x_2x_5m_3^4{+}99x_3x_5m_3^4{+}99x_4x_5m_3^4{+}54x_5^2m_3^4{-}72x_1m_1m_3^4\\
{-}72x_2m_1m_3^4{-}72x_3m_1m_3^4{-}72x_4m_1m_3^4{-}72x_5m_1m_3^4{+}21m_1^2m_3^4{-}72x_1m_2m_3^4\\
{-}72x_2m_2m_3^4{-}72x_3m_2m_3^4{-}72x_4m_2m_3^4{-}72x_5m_2m_3^4\\
{+}51m_1m_2m_3^4{+}21m_2^2m_3^4{-}1516x_1^5m_4{-}5404x_1^4x_2m_4{-}9094x_1^3x_2^2m_4{-}9094x_1^2x_2^3m_4\\
{-}5404x_1x_2^4m_4{-}1516x_2^5m_4{-}5404x_1^4x_3m_4{-}16178x_1^3x_2x_3m_4{-}22010x_1^2x_2^2x_3m_4\\
{-}16178x_1x_2^3x_3m_4{-}5404x_2^4x_3m_4{-}9094x_1^3x_3^2m_4{-}22010x_1^2x_2x_3^2m_4{-}22010x_1x_2^2x_3^2m_4\\
{-}9094x_2^3x_3^2m_4{-}9094x_1^2x_3^3m_4{-}16178x_1x_2x_3^3m_4{-}9094x_2^2x_3^3m_4{-}5404x_1x_3^4m_4\\
{-}5404x_2x_3^4m_4{-}1516x_3^5m_4{-}5404x_1^4x_4m_4{-}16178x_1^3x_2x_4m_4{-}22010x_1^2x_2^2x_4m_4\\
{-}16178x_1x_2^3x_4m_4{-}5404x_2^4x_4m_4{-}16178x_1^3x_3x_4m_4\\
{-}38946x_1^2x_2x_3x_4m_4{-}38946x_1x_2^2x_3x_4m_4{-}16178x_2^3x_3x_4m_4{-}22010x_1^2x_3^2x_4m_4\\
{-}38946x_1x_2x_3^2x_4m_4{-}22010x_2^2x_3^2x_4m_4{-}16178x_1x_3^3x_4m_4\\
{-}16178x_2x_3^3x_4m_4{-}5404x_3^4x_4m_4\\
{-}9094x_1^3x_4^2m_4{-}22010x_1^2x_2x_4^2m_4{-}22010x_1x_2^2x_4^2m_4{-}9094x_2^3x_4^2m_4\\
{-}22010x_1^2x_3x_4^2m_4{-}38946x_1x_2x_3x_4^2m_4{-}22010x_2^2x_3x_4^2m_4{-}22010x_1x_3^2x_4^2m_4\\
{-}22010x_2x_3^2x_4^2m_4{-}9094x_3^3x_4^2m_4\\
{-}9094x_1^2x_4^3m_4{-}16178x_1x_2x_4^3m_4{-}9094x_2^2x_4^3m_4{-}16178x_1x_3x_4^3m_4{-}16178x_2x_3x_4^3m_4\\
{-}9094x_3^2x_4^3m_4{-}5404x_1x_4^4m_4{-}5404x_2x_4^4m_4{-}5404x_3x_4^4m_4\\
{-}1516x_4^5m_4{-}5404x_1^4x_5m_4{-}16178x_1^3x_2x_5m_4\\
{-}22010x_1^2x_2^2x_5m_4{-}16178x_1x_2^3x_5m_4{-}5404x_2^4x_5m_4{-}16178x_1^3x_3x_5m_4\\
{-}38946x_1^2x_2x_3x_5m_4{-}38946x_1x_2^2x_3x_5m_4{-}16178x_2^3x_3x_5m_4{-}22010x_1^2x_3^2x_5m_4\\
{-}38946x_1x_2x_3^2x_5m_4{-}22010x_2^2x_3^2x_5m_4{-}16178x_1x_3^3x_5m_4{-}16178x_2x_3^3x_5m_4\\
{-}5404x_3^4x_5m_4{-}16178x_1^3x_4x_5m_4{-}38946x_1^2x_2x_4x_5m_4{-}38946x_1x_2^2x_4x_5m_4\\
{-}16178x_2^3x_4x_5m_4{-}38946x_1^2x_3x_4x_5m_4{-}68700x_1x_2x_3x_4x_5m_4{-}38946x_2^2x_3x_4x_5m_4\\
{-}38946x_1x_3^2x_4x_5m_4{-}38946x_2x_3^2x_4x_5m_4{-}16178x_3^3x_4x_5m_4{-}22010x_1^2x_4^2x_5m_4\\
{-}38946x_1x_2x_4^2x_5m_4{-}22010x_2^2x_4^2x_5m_4{-}38946x_1x_3x_4^2x_5m_4{-}38946x_2x_3x_4^2x_5m_4\\
{-}22010x_3^2x_4^2x_5m_4{-}16178x_1x_4^3x_5m_4{-}16178x_2x_4^3x_5m_4{-}16178x_3x_4^3x_5m_4\\
{-}5404x_4^4x_5m_4{-}9094x_1^3x_5^2m_4{-}22010x_1^2x_2x_5^2m_4\\
{-}22010x_1x_2^2x_5^2m_4{-}9094x_2^3x_5^2m_4{-}22010x_1^2x_3x_5^2m_4\\
{-}38946x_1x_2x_3x_5^2m_4{-}22010x_2^2x_3x_5^2m_4{-}22010x_1x_3^2x_5^2m_4{-}22010x_2x_3^2x_5^2m_4\\
{-}9094x_3^3x_5^2m_4{-}22010x_1^2x_4x_5^2m_4{-}38946x_1x_2x_4x_5^2m_4{-}22010x_2^2x_4x_5^2m_4\\
{-}38946x_1x_3x_4x_5^2m_4{-}38946x_2x_3x_4x_5^2m_4{-}22010x_3^2x_4x_5^2m_4{-}22010x_1x_4^2x_5^2m_4\\
{-}22010x_2x_4^2x_5^2m_4{-}22010x_3x_4^2x_5^2m_4{-}9094x_4^3x_5^2m_4{-}9094x_1^2x_5^3m_4\\
{-}16178x_1x_2x_5^3m_4{-}9094x_2^2x_5^3m_4{-}16178x_1x_3x_5^3m_4\\
{-}16178x_2x_3x_5^3m_4{-}9094x_3^2x_5^3m_4{-}16178x_1x_4x_5^3m_4\\
{-}16178x_2x_4x_5^3m_4{-}16178x_3x_4x_5^3m_4{-}9094x_4^2x_5^3m_4\\
{-}5404x_1x_5^4m_4{-}5404x_2x_5^4m_4{-}5404x_3x_5^4m_4{-}5404x_4x_5^4m_4{-}1516x_5^5m_4{+}3450x_1^4m_1m_4\\
{+}10664x_1^3x_2m_1m_4{+}14674x_1^2x_2^2m_1m_4{+}10664x_1x_2^3m_1m_4{+}3450x_2^4m_1m_4\\
{+}10664x_1^3x_3m_1m_4{+}26298x_1^2x_2x_3m_1m_4{+}26298x_1x_2^2x_3m_1m_4{+}10664x_2^3x_3m_1m_4\\
{+}14674x_1^2x_3^2m_1m_4{+}26298x_1x_2x_3^2m_1m_4{+}14674x_2^2x_3^2m_1m_4{+}10664x_1x_3^3m_1m_4\\
{+}10664x_2x_3^3m_1m_4{+}3450x_3^4m_1m_4{+}10664x_1^3x_4m_1m_4{+}26298x_1^2x_2x_4m_1m_4\\
{+}26298x_1x_2^2x_4m_1m_4{+}10664x_2^3x_4m_1m_4{+}26298x_1^2x_3x_4m_1m_4{+}47004x_1x_2x_3x_4m_1m_4\\
{+}26298x_2^2x_3x_4m_1m_4{+}26298x_1x_3^2x_4m_1m_4{+}26298x_2x_3^2x_4m_1m_4{+}10664x_3^3x_4m_1m_4\\
{+}14674x_1^2x_4^2m_1m_4{+}26298x_1x_2x_4^2m_1m_4{+}14674x_2^2x_4^2m_1m_4\\
{+}26298x_1x_3x_4^2m_1m_4{+}26298x_2x_3x_4^2m_1m_4{+}14674x_3^2x_4^2m_1m_4{+}10664x_1x_4^3m_1m_4\\
{+}10664x_2x_4^3m_1m_4{+}10664x_3x_4^3m_1m_4{+}3450x_4^4m_1m_4{+}10664x_1^3x_5m_1m_4\\
{+}26298x_1^2x_2x_5m_1m_4{+}26298x_1x_2^2x_5m_1m_4{+}10664x_2^3x_5m_1m_4{+}26298x_1^2x_3x_5m_1m_4\\
{+}47004x_1x_2x_3x_5m_1m_4{+}26298x_2^2x_3x_5m_1m_4{+}26298x_1x_3^2x_5m_1m_4{+}26298x_2x_3^2x_5m_1m_4\\
{+}10664x_3^3x_5m_1m_4{+}26298x_1^2x_4x_5m_1m_4\\
{+}47004x_1x_2x_4x_5m_1m_4{+}26298x_2^2x_4x_5m_1m_4{+}47004x_1x_3x_4x_5m_1m_4{+}47004x_2x_3x_4x_5m_1m_4\\
{+}26298x_3^2x_4x_5m_1m_4{+}26298x_1x_4^2x_5m_1m_4{+}26298x_2x_4^2x_5m_1m_4{+}26298x_3x_4^2x_5m_1m_4\\
{+}10664x_4^3x_5m_1m_4{+}14674x_1^2x_5^2m_1m_4{+}26298x_1x_2x_5^2m_1m_4{+}14674x_2^2x_5^2m_1m_4\\
{+}26298x_1x_3x_5^2m_1m_4{+}26298x_2x_3x_5^2m_1m_4{+}14674x_3^2x_5^2m_1m_4{+}26298x_1x_4x_5^2m_1m_4\\
{+}26298x_2x_4x_5^2m_1m_4{+}26298x_3x_4x_5^2m_1m_4\\
{+}14674x_4^2x_5^2m_1m_4{+}10664x_1x_5^3m_1m_4{+}10664x_2x_5^3m_1m_4{+}10664x_3x_5^3m_1m_4\\
{+}10664x_4x_5^3m_1m_4{+}3450x_5^4m_1m_4{-}2663x_1^3m_1^2m_4\\
{-}6694x_1^2x_2m_1^2m_4{-}6694x_1x_2^2m_1^2m_4\\
{-}2663x_2^3m_1^2m_4{-}6694x_1^2x_3m_1^2m_4{-}12093x_1x_2x_3m_1^2m_4{-}6694x_2^2x_3m_1^2m_4\\
{-}6694x_1x_3^2m_1^2m_4{-}6694x_2x_3^2m_1^2m_4{-}2663x_3^3m_1^2m_4{-}6694x_1^2x_4m_1^2m_4\\
{-}12093x_1x_2x_4m_1^2m_4{-}6694x_2^2x_4m_1^2m_4{-}12093x_1x_3x_4m_1^2m_4{-}12093x_2x_3x_4m_1^2m_4\\
{-}6694x_3^2x_4m_1^2m_4{-}6694x_1x_4^2m_1^2m_4{-}6694x_2x_4^2m_1^2m_4{-}6694x_3x_4^2m_1^2m_4\\
{-}2663x_4^3m_1^2m_4{-}6694x_1^2x_5m_1^2m_4{-}12093x_1x_2x_5m_1^2m_4\\
{-}6694x_2^2x_5m_1^2m_4{-}12093x_1x_3x_5m_1^2m_4{-}12093x_2x_3x_5m_1^2m_4{-}6694x_3^2x_5m_1^2m_4\\
{-}12093x_1x_4x_5m_1^2m_4{-}12093x_2x_4x_5m_1^2m_4{-}12093x_3x_4x_5m_1^2m_4{-}6694x_4^2x_5m_1^2m_4\\
{-}6694x_1x_5^2m_1^2m_4{-}6694x_2x_5^2m_1^2m_4{-}6694x_3x_5^2m_1^2m_4{-}6694x_4x_5^2m_1^2m_4\\
{-}2663x_5^3m_1^2m_4{+}798x_1^2m_1^3m_4{+}1458x_1x_2m_1^3m_4{+}798x_2^2m_1^3m_4{+}1458x_1x_3m_1^3m_4\\
{+}1458x_2x_3m_1^3m_4{+}798x_3^2m_1^3m_4{+}1458x_1x_4m_1^3m_4{+}1458x_2x_4m_1^3m_4{+}1458x_3x_4m_1^3m_4\\
{+}798x_4^2m_1^3m_4{+}1458x_1x_5m_1^3m_4{+}1458x_2x_5m_1^3m_4{+}1458x_3x_5m_1^3m_4{+}1458x_4x_5m_1^3m_4\\
{+}798x_5^2m_1^3m_4{-}72x_1m_1^4m_4{-}72x_2m_1^4m_4{-}72x_3m_1^4m_4{-}72x_4m_1^4m_4{-}72x_5m_1^4m_4\\
{+}3450x_1^4m_2m_4{+}10664x_1^3x_2m_2m_4{+}14674x_1^2x_2^2m_2m_4{+}10664x_1x_2^3m_2m_4{+}3450x_2^4m_2m_4\\
{+}10664x_1^3x_3m_2m_4{+}26298x_1^2x_2x_3m_2m_4{+}26298x_1x_2^2x_3m_2m_4{+}10664x_2^3x_3m_2m_4\\
{+}14674x_1^2x_3^2m_2m_4{+}26298x_1x_2x_3^2m_2m_4{+}14674x_2^2x_3^2m_2m_4\\
{+}10664x_1x_3^3m_2m_4{+}10664x_2x_3^3m_2m_4\\
{+}3450x_3^4m_2m_4{+}10664x_1^3x_4m_2m_4{+}26298x_1^2x_2x_4m_2m_4{+}26298x_1x_2^2x_4m_2m_4\\
{+}10664x_2^3x_4m_2m_4{+}26298x_1^2x_3x_4m_2m_4{+}47004x_1x_2x_3x_4m_2m_4\\
{+}26298x_2^2x_3x_4m_2m_4{+}26298x_1x_3^2x_4m_2m_4{+}26298x_2x_3^2x_4m_2m_4\\
{+}10664x_3^3x_4m_2m_4{+}14674x_1^2x_4^2m_2m_4\\
{+}26298x_1x_2x_4^2m_2m_4{+}14674x_2^2x_4^2m_2m_4{+}26298x_1x_3x_4^2m_2m_4{+}26298x_2x_3x_4^2m_2m_4\\
{+}14674x_3^2x_4^2m_2m_4{+}10664x_1x_4^3m_2m_4{+}10664x_2x_4^3m_2m_4{+}10664x_3x_4^3m_2m_4\\
{+}3450x_4^4m_2m_4{+}10664x_1^3x_5m_2m_4{+}26298x_1^2x_2x_5m_2m_4{+}26298x_1x_2^2x_5m_2m_4\\
{+}10664x_2^3x_5m_2m_4{+}26298x_1^2x_3x_5m_2m_4{+}47004x_1x_2x_3x_5m_2m_4{+}26298x_2^2x_3x_5m_2m_4\\
{+}26298x_1x_3^2x_5m_2m_4{+}26298x_2x_3^2x_5m_2m_4{+}10664x_3^3x_5m_2m_4\\
{+}26298x_1^2x_4x_5m_2m_4{+}47004x_1x_2x_4x_5m_2m_4\\
{+}26298x_2^2x_4x_5m_2m_4{+}47004x_1x_3x_4x_5m_2m_4{+}47004x_2x_3x_4x_5m_2m_4{+}26298x_3^2x_4x_5m_2m_4\\
{+}26298x_1x_4^2x_5m_2m_4{+}26298x_2x_4^2x_5m_2m_4{+}26298x_3x_4^2x_5m_2m_4{+}10664x_4^3x_5m_2m_4\\
{+}14674x_1^2x_5^2m_2m_4{+}26298x_1x_2x_5^2m_2m_4\\
{+}14674x_2^2x_5^2m_2m_4{+}26298x_1x_3x_5^2m_2m_4{+}26298x_2x_3x_5^2m_2m_4{+}14674x_3^2x_5^2m_2m_4\\
{+}26298x_1x_4x_5^2m_2m_4{+}26298x_2x_4x_5^2m_2m_4{+}26298x_3x_4x_5^2m_2m_4{+}14674x_4^2x_5^2m_2m_4\\
{+}10664x_1x_5^3m_2m_4{+}10664x_2x_5^3m_2m_4{+}10664x_3x_5^3m_2m_4{+}10664x_4x_5^3m_2m_4{+}3450x_5^4m_2m_4\\
{-}6912x_1^3m_1m_2m_4{-}17390x_1^2x_2m_1m_2m_4{-}17390x_1x_2^2m_1m_2m_4{-}6912x_2^3m_1m_2m_4\\
{-}17390x_1^2x_3m_1m_2m_4{-}31434x_1x_2x_3m_1m_2m_4{-}17390x_2^2x_3m_1m_2m_4{-}17390x_1x_3^2m_1m_2m_4\\
{-}17390x_2x_3^2m_1m_2m_4{-}6912x_3^3m_1m_2m_4{-}17390x_1^2x_4m_1m_2m_4{-}31434x_1x_2x_4m_1m_2m_4\\
{-}17390x_2^2x_4m_1m_2m_4{-}31434x_1x_3x_4m_1m_2m_4{-}31434x_2x_3x_4m_1m_2m_4{-}17390x_3^2x_4m_1m_2m_4\\
{-}17390x_1x_4^2m_1m_2m_4{-}17390x_2x_4^2m_1m_2m_4{-}17390x_3x_4^2m_1m_2m_4{-}6912x_4^3m_1m_2m_4\\
{-}17390x_1^2x_5m_1m_2m_4{-}31434x_1x_2x_5m_1m_2m_4{-}17390x_2^2x_5m_1m_2m_4\\
{-}31434x_1x_3x_5m_1m_2m_4{-}31434x_2x_3x_5m_1m_2m_4\\
{-}17390x_3^2x_5m_1m_2m_4{-}31434x_1x_4x_5m_1m_2m_4\\
{-}31434x_2x_4x_5m_1m_2m_4{-}31434x_3x_4x_5m_1m_2m_4{-}17390x_4^2x_5m_1m_2m_4{-}17390x_1x_5^2m_1m_2m_4\\
{-}17390x_2x_5^2m_1m_2m_4{-}17390x_3x_5^2m_1m_2m_4{-}17390x_4x_5^2m_1m_2m_4{-}6912x_5^3m_1m_2m_4\\
{+}4452x_1^2m_1^2m_2m_4{+}8122x_1x_2m_1^2m_2m_4{+}4452x_2^2m_1^2m_2m_4{+}8122x_1x_3m_1^2m_2m_4\\
{+}8122x_2x_3m_1^2m_2m_4{+}4452x_3^2m_1^2m_2m_4{+}8122x_1x_4m_1^2m_2m_4{+}8122x_2x_4m_1^2m_2m_4\\
{+}8122x_3x_4m_1^2m_2m_4{+}4452x_4^2m_1^2m_2m_4{+}8122x_1x_5m_1^2m_2m_4{+}8122x_2x_5m_1^2m_2m_4\\
{+}8122x_3x_5m_1^2m_2m_4{+}8122x_4x_5m_1^2m_2m_4{+}4452x_5^2m_1^2m_2m_4{-}1002x_1m_1^3m_2m_4\\
{-}1002x_2m_1^3m_2m_4{-}1002x_3m_1^3m_2m_4{-}1002x_4m_1^3m_2m_4{-}1002x_5m_1^3m_2m_4{+}51m_1^4m_2m_4\\
{-}2663x_1^3m_2^2m_4{-}6694x_1^2x_2m_2^2m_4{-}6694x_1x_2^2m_2^2m_4{-}2663x_2^3m_2^2m_4\\
{-}6694x_1^2x_3m_2^2m_4{-}12093x_1x_2x_3m_2^2m_4{-}6694x_2^2x_3m_2^2m_4{-}6694x_1x_3^2m_2^2m_4\\
{-}6694x_2x_3^2m_2^2m_4{-}2663x_3^3m_2^2m_4{-}6694x_1^2x_4m_2^2m_4{-}12093x_1x_2x_4m_2^2m_4\\
{-}6694x_2^2x_4m_2^2m_4{-}12093x_1x_3x_4m_2^2m_4{-}12093x_2x_3x_4m_2^2m_4\\
{-}6694x_3^2x_4m_2^2m_4{-}6694x_1x_4^2m_2^2m_4{-}6694x_2x_4^2m_2^2m_4{-}6694x_3x_4^2m_2^2m_4\\
{-}2663x_4^3m_2^2m_4{-}6694x_1^2x_5m_2^2m_4{-}12093x_1x_2x_5m_2^2m_4\\
{-}6694x_2^2x_5m_2^2m_4{-}12093x_1x_3x_5m_2^2m_4{-}12093x_2x_3x_5m_2^2m_4{-}6694x_3^2x_5m_2^2m_4\\
{-}12093x_1x_4x_5m_2^2m_4{-}12093x_2x_4x_5m_2^2m_4{-}12093x_3x_4x_5m_2^2m_4{-}6694x_4^2x_5m_2^2m_4\\
{-}6694x_1x_5^2m_2^2m_4{-}6694x_2x_5^2m_2^2m_4{-}6694x_3x_5^2m_2^2m_4{-}6694x_4x_5^2m_2^2m_4\\
{-}2663x_5^3m_2^2m_4{+}4452x_1^2m_1m_2^2m_4{+}8122x_1x_2m_1m_2^2m_4{+}4452x_2^2m_1m_2^2m_4\\
{+}8122x_1x_3m_1m_2^2m_4{+}8122x_2x_3m_1m_2^2m_4{+}4452x_3^2m_1m_2^2m_4{+}8122x_1x_4m_1m_2^2m_4\\
{+}8122x_2x_4m_1m_2^2m_4{+}8122x_3x_4m_1m_2^2m_4{+}4452x_4^2m_1m_2^2m_4{+}8122x_1x_5m_1m_2^2m_4\\
{+}8122x_2x_5m_1m_2^2m_4{+}8122x_3x_5m_1m_2^2m_4{+}8122x_4x_5m_1m_2^2m_4\\
{+}4452x_5^2m_1m_2^2m_4{-}2158x_1m_1^2m_2^2m_4{-}2158x_2m_1^2m_2^2m_4{-}2158x_3m_1^2m_2^2m_4\\
{-}2158x_4m_1^2m_2^2m_4{-}2158x_5m_1^2m_2^2m_4{+}279m_1^3m_2^2m_4{+}798x_1^2m_2^3m_4\\
{+}1458x_1x_2m_2^3m_4{+}798x_2^2m_2^3m_4{+}1458x_1x_3m_2^3m_4{+}1458x_2x_3m_2^3m_4\\
{+}798x_3^2m_2^3m_4{+}1458x_1x_4m_2^3m_4{+}1458x_2x_4m_2^3m_4{+}1458x_3x_4m_2^3m_4{+}798x_4^2m_2^3m_4\\
{+}1458x_1x_5m_2^3m_4{+}1458x_2x_5m_2^3m_4{+}1458x_3x_5m_2^3m_4{+}1458x_4x_5m_2^3m_4{+}798x_5^2m_2^3m_4\\
{-}1002x_1m_1m_2^3m_4{-}1002x_2m_1m_2^3m_4{-}1002x_3m_1m_2^3m_4{-}1002x_4m_1m_2^3m_4\\
{-}1002x_5m_1m_2^3m_4{+}279m_1^2m_2^3m_4{-}72x_1m_2^4m_4{-}72x_2m_2^4m_4{-}72x_3m_2^4m_4{-}72x_4m_2^4m_4\\
{-}72x_5m_2^4m_4{+}51m_1m_2^4m_4{+}3450x_1^4m_3m_4{+}10664x_1^3x_2m_3m_4{+}14674x_1^2x_2^2m_3m_4\\
{+}10664x_1x_2^3m_3m_4{+}3450x_2^4m_3m_4{+}10664x_1^3x_3m_3m_4{+}26298x_1^2x_2x_3m_3m_4\\
{+}26298x_1x_2^2x_3m_3m_4{+}10664x_2^3x_3m_3m_4{+}14674x_1^2x_3^2m_3m_4\\
{+}26298x_1x_2x_3^2m_3m_4{+}14674x_2^2x_3^2m_3m_4\\
{+}10664x_1x_3^3m_3m_4{+}10664x_2x_3^3m_3m_4{+}3450x_3^4m_3m_4{+}10664x_1^3x_4m_3m_4\\
{+}26298x_1^2x_2x_4m_3m_4{+}26298x_1x_2^2x_4m_3m_4{+}10664x_2^3x_4m_3m_4{+}26298x_1^2x_3x_4m_3m_4\\
{+}47004x_1x_2x_3x_4m_3m_4{+}26298x_2^2x_3x_4m_3m_4{+}26298x_1x_3^2x_4m_3m_4{+}26298x_2x_3^2x_4m_3m_4\\
{+}10664x_3^3x_4m_3m_4{+}14674x_1^2x_4^2m_3m_4{+}26298x_1x_2x_4^2m_3m_4{+}14674x_2^2x_4^2m_3m_4\\
{+}26298x_1x_3x_4^2m_3m_4{+}26298x_2x_3x_4^2m_3m_4{+}14674x_3^2x_4^2m_3m_4{+}10664x_1x_4^3m_3m_4\\
{+}10664x_2x_4^3m_3m_4{+}10664x_3x_4^3m_3m_4{+}3450x_4^4m_3m_4\\
{+}10664x_1^3x_5m_3m_4{+}26298x_1^2x_2x_5m_3m_4\\
{+}26298x_1x_2^2x_5m_3m_4{+}10664x_2^3x_5m_3m_4{+}26298x_1^2x_3x_5m_3m_4{+}47004x_1x_2x_3x_5m_3m_4\\
{+}26298x_2^2x_3x_5m_3m_4{+}26298x_1x_3^2x_5m_3m_4{+}26298x_2x_3^2x_5m_3m_4{+}10664x_3^3x_5m_3m_4\\{+}26298x_1^2x_4x_5m_3m_4{+}47004x_1x_2x_4x_5m_3m_4{+}26298x_2^2x_4x_5m_3m_4{+}47004x_1x_3x_4x_5m_3m_4\\{+}47004x_2x_3x_4x_5m_3m_4{+}26298x_3^2x_4x_5m_3m_4{+}26298x_1x_4^2x_5m_3m_4{+}26298x_2x_4^2x_5m_3m_4\\{+}26298x_3x_4^2x_5m_3m_4{+}10664x_4^3x_5m_3m_4{+}14674x_1^2x_5^2m_3m_4{+}26298x_1x_2x_5^2m_3m_4\\
{+}14674x_2^2x_5^2m_3m_4{+}26298x_1x_3x_5^2m_3m_4{+}26298x_2x_3x_5^2m_3m_4\\
{+}14674x_3^2x_5^2m_3m_4{+}26298x_1x_4x_5^2m_3m_4\\
{+}26298x_2x_4x_5^2m_3m_4{+}26298x_3x_4x_5^2m_3m_4{+}14674x_4^2x_5^2m_3m_4{+}10664x_1x_5^3m_3m_4\\
{+}10664x_2x_5^3m_3m_4{+}10664x_3x_5^3m_3m_4{+}10664x_4x_5^3m_3m_4{+}3450x_5^4m_3m_4{-}6912x_1^3m_1m_3m_4\\
{-}17390x_1^2x_2m_1m_3m_4{-}17390x_1x_2^2m_1m_3m_4{-}6912x_2^3m_1m_3m_4{-}17390x_1^2x_3m_1m_3m_4\\
{-}31434x_1x_2x_3m_1m_3m_4{-}17390x_2^2x_3m_1m_3m_4{-}17390x_1x_3^2m_1m_3m_4{-}17390x_2x_3^2m_1m_3m_4\\{-}6912x_3^3m_1m_3m_4{-}17390x_1^2x_4m_1m_3m_4{-}31434x_1x_2x_4m_1m_3m_4{-}17390x_2^2x_4m_1m_3m_4\\{-}31434x_1x_3x_4m_1m_3m_4{-}31434x_2x_3x_4m_1m_3m_4{-}17390x_3^2x_4m_1m_3m_4{-}17390x_1x_4^2m_1m_3m_4\\{-}17390x_2x_4^2m_1m_3m_4{-}17390x_3x_4^2m_1m_3m_4{-}6912x_4^3m_1m_3m_4{-}17390x_1^2x_5m_1m_3m_4\\
{-}31434x_1x_2x_5m_1m_3m_4{-}17390x_2^2x_5m_1m_3m_4\\
{-}31434x_1x_3x_5m_1m_3m_4{-}31434x_2x_3x_5m_1m_3m_4\\
{-}17390x_3^2x_5m_1m_3m_4{-}31434x_1x_4x_5m_1m_3m_4\\
{-}31434x_2x_4x_5m_1m_3m_4{-}31434x_3x_4x_5m_1m_3m_4{-}17390x_4^2x_5m_1m_3m_4{-}17390x_1x_5^2m_1m_3m_4\\
{-}17390x_2x_5^2m_1m_3m_4{-}17390x_3x_5^2m_1m_3m_4{-}17390x_4x_5^2m_1m_3m_4{-}6912x_5^3m_1m_3m_4\\
{+}4452x_1^2m_1^2m_3m_4{+}8122x_1x_2m_1^2m_3m_4{+}4452x_2^2m_1^2m_3m_4{+}8122x_1x_3m_1^2m_3m_4\\
{+}8122x_2x_3m_1^2m_3m_4{+}4452x_3^2m_1^2m_3m_4{+}8122x_1x_4m_1^2m_3m_4{+}8122x_2x_4m_1^2m_3m_4\\
{+}8122x_3x_4m_1^2m_3m_4{+}4452x_4^2m_1^2m_3m_4{+}8122x_1x_5m_1^2m_3m_4{+}8122x_2x_5m_1^2m_3m_4\\
{+}8122x_3x_5m_1^2m_3m_4{+}8122x_4x_5m_1^2m_3m_4{+}4452x_5^2m_1^2m_3m_4{-}1002x_1m_1^3m_3m_4\\
{-}1002x_2m_1^3m_3m_4{-}1002x_3m_1^3m_3m_4{-}1002x_4m_1^3m_3m_4{-}1002x_5m_1^3m_3m_4{+}51m_1^4m_3m_4\\
{-}6912x_1^3m_2m_3m_4{-}17390x_1^2x_2m_2m_3m_4{-}17390x_1x_2^2m_2m_3m_4\\
{-}6912x_2^3m_2m_3m_4{-}17390x_1^2x_3m_2m_3m_4{-}31434x_1x_2x_3m_2m_3m_4{-}17390x_2^2x_3m_2m_3m_4\\
{-}17390x_1x_3^2m_2m_3m_4{-}17390x_2x_3^2m_2m_3m_4{-}6912x_3^3m_2m_3m_4{-}17390x_1^2x_4m_2m_3m_4\\
{-}31434x_1x_2x_4m_2m_3m_4{-}17390x_2^2x_4m_2m_3m_4{-}31434x_1x_3x_4m_2m_3m_4\\
{-}31434x_2x_3x_4m_2m_3m_4{-}17390x_3^2x_4m_2m_3m_4{-}17390x_1x_4^2m_2m_3m_4\\
{-}17390x_2x_4^2m_2m_3m_4{-}17390x_3x_4^2m_2m_3m_4{-}6912x_4^3m_2m_3m_4\\
{-}17390x_1^2x_5m_2m_3m_4{-}31434x_1x_2x_5m_2m_3m_4{-}17390x_2^2x_5m_2m_3m_4\\
{-}31434x_1x_3x_5m_2m_3m_4{-}31434x_2x_3x_5m_2m_3m_4{-}17390x_3^2x_5m_2m_3m_4\\
{-}31434x_1x_4x_5m_2m_3m_4{-}31434x_2x_4x_5m_2m_3m_4{-}31434x_3x_4x_5m_2m_3m_4\\
{-}17390x_4^2x_5m_2m_3m_4{-}17390x_1x_5^2m_2m_3m_4{-}17390x_2x_5^2m_2m_3m_4\\
{-}17390x_3x_5^2m_2m_3m_4{-}17390x_4x_5^2m_2m_3m_4{-}6912x_5^3m_2m_3m_4{+}11336x_1^2m_1m_2m_3m_4\\
{+}20672x_1x_2m_1m_2m_3m_4{+}11336x_2^2m_1m_2m_3m_4\\
{+}20672x_1x_3m_1m_2m_3m_4{+}20672x_2x_3m_1m_2m_3m_4\\
{+}11336x_3^2m_1m_2m_3m_4{+}20672x_1x_4m_1m_2m_3m_4\\
{+}20672x_2x_4m_1m_2m_3m_4{+}20672x_3x_4m_1m_2m_3m_4{+}11336x_4^2m_1m_2m_3m_4\\
{+}20672x_1x_5m_1m_2m_3m_4{+}20672x_2x_5m_1m_2m_3m_4{+}20672x_3x_5m_1m_2m_3m_4\\
{+}20672x_4x_5m_1m_2m_3m_4{+}11336x_5^2m_1m_2m_3m_4\\
{-}5394x_1m_1^2m_2m_3m_4{-}5394x_2m_1^2m_2m_3m_4\\
{-}5394x_3m_1^2m_2m_3m_4{-}5394x_4m_1^2m_2m_3m_4{-}5394x_5m_1^2m_2m_3m_4{+}684m_1^3m_2m_3m_4\\
{+}4452x_1^2m_2^2m_3m_4{+}8122x_1x_2m_2^2m_3m_4{+}4452x_2^2m_2^2m_3m_4\\
{+}8122x_1x_3m_2^2m_3m_4{+}8122x_2x_3m_2^2m_3m_4\\
{+}4452x_3^2m_2^2m_3m_4{+}8122x_1x_4m_2^2m_3m_4{+}8122x_2x_4m_2^2m_3m_4{+}8122x_3x_4m_2^2m_3m_4\\
{+}4452x_4^2m_2^2m_3m_4{+}8122x_1x_5m_2^2m_3m_4{+}8122x_2x_5m_2^2m_3m_4{+}8122x_3x_5m_2^2m_3m_4\\
{+}8122x_4x_5m_2^2m_3m_4{+}4452x_5^2m_2^2m_3m_4{-}5394x_1m_1m_2^2m_3m_4{-}5394x_2m_1m_2^2m_3m_4\\
{-}5394x_3m_1m_2^2m_3m_4{-}5394x_4m_1m_2^2m_3m_4{-}5394x_5m_1m_2^2m_3m_4{+}1456m_1^2m_2^2m_3m_4\\
{-}1002x_1m_2^3m_3m_4{-}1002x_2m_2^3m_3m_4{-}1002x_3m_2^3m_3m_4{-}1002x_4m_2^3m_3m_4\\
{-}1002x_5m_2^3m_3m_4{+}684m_1m_2^3m_3m_4{+}51m_2^4m_3m_4{-}2663x_1^3m_3^2m_4{-}6694x_1^2x_2m_3^2m_4\\
{-}6694x_1x_2^2m_3^2m_4{-}2663x_2^3m_3^2m_4{-}6694x_1^2x_3m_3^2m_4{-}12093x_1x_2x_3m_3^2m_4\\
{-}6694x_2^2x_3m_3^2m_4{-}6694x_1x_3^2m_3^2m_4{-}6694x_2x_3^2m_3^2m_4{-}2663x_3^3m_3^2m_4\\
{-}6694x_1^2x_4m_3^2m_4{-}12093x_1x_2x_4m_3^2m_4\\
{-}6694x_2^2x_4m_3^2m_4{-}12093x_1x_3x_4m_3^2m_4{-}12093x_2x_3x_4m_3^2m_4{-}6694x_3^2x_4m_3^2m_4\\
{-}6694x_1x_4^2m_3^2m_4{-}6694x_2x_4^2m_3^2m_4{-}6694x_3x_4^2m_3^2m_4{-}2663x_4^3m_3^2m_4\\
{-}6694x_1^2x_5m_3^2m_4{-}12093x_1x_2x_5m_3^2m_4{-}6694x_2^2x_5m_3^2m_4{-}12093x_1x_3x_5m_3^2m_4\\
{-}12093x_2x_3x_5m_3^2m_4{-}6694x_3^2x_5m_3^2m_4{-}12093x_1x_4x_5m_3^2m_4{-}12093x_2x_4x_5m_3^2m_4\\
{-}12093x_3x_4x_5m_3^2m_4{-}6694x_4^2x_5m_3^2m_4{-}6694x_1x_5^2m_3^2m_4{-}6694x_2x_5^2m_3^2m_4\\
{-}6694x_3x_5^2m_3^2m_4{-}6694x_4x_5^2m_3^2m_4{-}2663x_5^3m_3^2m_4{+}4452x_1^2m_1m_3^2m_4\\
{+}8122x_1x_2m_1m_3^2m_4{+}4452x_2^2m_1m_3^2m_4\\
{+}8122x_1x_3m_1m_3^2m_4{+}8122x_2x_3m_1m_3^2m_4{+}4452x_3^2m_1m_3^2m_4{+}8122x_1x_4m_1m_3^2m_4\\
{+}8122x_2x_4m_1m_3^2m_4{+}8122x_3x_4m_1m_3^2m_4{+}4452x_4^2m_1m_3^2m_4{+}8122x_1x_5m_1m_3^2m_4\\
{+}8122x_2x_5m_1m_3^2m_4{+}8122x_3x_5m_1m_3^2m_4{+}8122x_4x_5m_1m_3^2m_4\\
{+}4452x_5^2m_1m_3^2m_4{-}2158x_1m_1^2m_3^2m_4{-}2158x_2m_1^2m_3^2m_4{-}2158x_3m_1^2m_3^2m_4\\
{-}2158x_4m_1^2m_3^2m_4{-}2158x_5m_1^2m_3^2m_4{+}279m_1^3m_3^2m_4{+}4452x_1^2m_2m_3^2m_4\\
{+}8122x_1x_2m_2m_3^2m_4{+}4452x_2^2m_2m_3^2m_4{+}8122x_1x_3m_2m_3^2m_4\\
{+}8122x_2x_3m_2m_3^2m_4{+}4452x_3^2m_2m_3^2m_4\\
{+}8122x_1x_4m_2m_3^2m_4{+}8122x_2x_4m_2m_3^2m_4{+}8122x_3x_4m_2m_3^2m_4{+}4452x_4^2m_2m_3^2m_4\\
{+}8122x_1x_5m_2m_3^2m_4{+}8122x_2x_5m_2m_3^2m_4{+}8122x_3x_5m_2m_3^2m_4{+}8122x_4x_5m_2m_3^2m_4\\
{+}4452x_5^2m_2m_3^2m_4{-}5394x_1m_1m_2m_3^2m_4{-}5394x_2m_1m_2m_3^2m_4{-}5394x_3m_1m_2m_3^2m_4\\{-}5394x_4m_1m_2m_3^2m_4{-}5394x_5m_1m_2m_3^2m_4{+}1456m_1^2m_2m_3^2m_4\\
{-}2158x_1m_2^2m_3^2m_4{-}2158x_2m_2^2m_3^2m_4\\
{-}2158x_3m_2^2m_3^2m_4{-}2158x_4m_2^2m_3^2m_4{-}2158x_5m_2^2m_3^2m_4{+}1456m_1m_2^2m_3^2m_4\\
{+}279m_2^3m_3^2m_4{+}798x_1^2m_3^3m_4{+}1458x_1x_2m_3^3m_4{+}798x_2^2m_3^3m_4{+}1458x_1x_3m_3^3m_4\\
{+}1458x_2x_3m_3^3m_4{+}798x_3^2m_3^3m_4{+}1458x_1x_4m_3^3m_4{+}1458x_2x_4m_3^3m_4{+}1458x_3x_4m_3^3m_4\\
{+}798x_4^2m_3^3m_4{+}1458x_1x_5m_3^3m_4{+}1458x_2x_5m_3^3m_4{+}1458x_3x_5m_3^3m_4{+}1458x_4x_5m_3^3m_4\\
{+}798x_5^2m_3^3m_4{-}1002x_1m_1m_3^3m_4{-}1002x_2m_1m_3^3m_4{-}1002x_3m_1m_3^3m_4{-}1002x_4m_1m_3^3m_4\\
{-}1002x_5m_1m_3^3m_4{+}279m_1^2m_3^3m_4{-}1002x_1m_2m_3^3m_4{-}1002x_2m_2m_3^3m_4{-}1002x_3m_2m_3^3m_4\\
{-}1002x_4m_2m_3^3m_4{-}1002x_5m_2m_3^3m_4\\
{+}684m_1m_2m_3^3m_4{+}279m_2^2m_3^3m_4{-}72x_1m_3^4m_4{-}72x_2m_3^4m_4{-}72x_3m_3^4m_4{-}72x_4m_3^4m_4\\
{-}72x_5m_3^4m_4{+}51m_1m_3^4m_4{+}51m_2m_3^4m_4{+}1305x_1^4m_4^2{+}4020x_1^3x_2m_4^2\\
{+}5523x_1^2x_2^2m_4^2{+}4020x_1x_2^3m_4^2{+}1305x_2^4m_4^2{+}4020x_1^3x_3m_4^2{+}9883x_1^2x_2x_3m_4^2\\
{+}9883x_1x_2^2x_3m_4^2{+}4020x_2^3x_3m_4^2{+}5523x_1^2x_3^2m_4^2{+}9883x_1x_2x_3^2m_4^2{+}5523x_2^2x_3^2m_4^2\\
{+}4020x_1x_3^3m_4^2{+}4020x_2x_3^3m_4^2{+}1305x_3^4m_4^2{+}4020x_1^3x_4m_4^2{+}9883x_1^2x_2x_4m_4^2\\
{+}9883x_1x_2^2x_4m_4^2{+}4020x_2^3x_4m_4^2{+}9883x_1^2x_3x_4m_4^2{+}17634x_1x_2x_3x_4m_4^2\\
{+}9883x_2^2x_3x_4m_4^2{+}9883x_1x_3^2x_4m_4^2{+}9883x_2x_3^2x_4m_4^2{+}4020x_3^3x_4m_4^2\\
{+}5523x_1^2x_4^2m_4^2{+}9883x_1x_2x_4^2m_4^2{+}5523x_2^2x_4^2m_4^2{+}9883x_1x_3x_4^2m_4^2\\
{+}9883x_2x_3x_4^2m_4^2{+}5523x_3^2x_4^2m_4^2{+}4020x_1x_4^3m_4^2{+}4020x_2x_4^3m_4^2{+}4020x_3x_4^3m_4^2\\
{+}1305x_4^4m_4^2{+}4020x_1^3x_5m_4^2{+}9883x_1^2x_2x_5m_4^2{+}9883x_1x_2^2x_5m_4^2{+}4020x_2^3x_5m_4^2\\
{+}9883x_1^2x_3x_5m_4^2{+}17634x_1x_2x_3x_5m_4^2{+}9883x_2^2x_3x_5m_4^2{+}9883x_1x_3^2x_5m_4^2\\
{+}9883x_2x_3^2x_5m_4^2{+}4020x_3^3x_5m_4^2{+}9883x_1^2x_4x_5m_4^2{+}17634x_1x_2x_4x_5m_4^2\\
{+}9883x_2^2x_4x_5m_4^2{+}17634x_1x_3x_4x_5m_4^2{+}17634x_2x_3x_4x_5m_4^2{+}9883x_3^2x_4x_5m_4^2\\
{+}9883x_1x_4^2x_5m_4^2{+}9883x_2x_4^2x_5m_4^2{+}9883x_3x_4^2x_5m_4^2\\
{+}4020x_4^3x_5m_4^2{+}5523x_1^2x_5^2m_4^2{+}9883x_1x_2x_5^2m_4^2{+}5523x_2^2x_5^2m_4^2\\
{+}9883x_1x_3x_5^2m_4^2{+}9883x_2x_3x_5^2m_4^2{+}5523x_3^2x_5^2m_4^2{+}9883x_1x_4x_5^2m_4^2\\
{+}9883x_2x_4x_5^2m_4^2{+}9883x_3x_4x_5^2m_4^2{+}5523x_4^2x_5^2m_4^2\\
{+}4020x_1x_5^3m_4^2{+}4020x_2x_5^3m_4^2\\
{+}4020x_3x_5^3m_4^2{+}4020x_4x_5^3m_4^2{+}1305x_5^4m_4^2{-}2663x_1^3m_1m_4^2{-}6694x_1^2x_2m_1m_4^2\\
{-}6694x_1x_2^2m_1m_4^2{-}2663x_2^3m_1m_4^2{-}6694x_1^2x_3m_1m_4^2{-}12093x_1x_2x_3m_1m_4^2\\
{-}6694x_2^2x_3m_1m_4^2{-}6694x_1x_3^2m_1m_4^2{-}6694x_2x_3^2m_1m_4^2{-}2663x_3^3m_1m_4^2\\
{-}6694x_1^2x_4m_1m_4^2{-}12093x_1x_2x_4m_1m_4^2{-}6694x_2^2x_4m_1m_4^2{-}12093x_1x_3x_4m_1m_4^2\\
{-}12093x_2x_3x_4m_1m_4^2{-}6694x_3^2x_4m_1m_4^2\\
{-}6694x_1x_4^2m_1m_4^2{-}6694x_2x_4^2m_1m_4^2{-}6694x_3x_4^2m_1m_4^2\\
{-}2663x_4^3m_1m_4^2{-}6694x_1^2x_5m_1m_4^2{-}12093x_1x_2x_5m_1m_4^2\\
{-}6694x_2^2x_5m_1m_4^2{-}12093x_1x_3x_5m_1m_4^2{-}12093x_2x_3x_5m_1m_4^2\\
{-}6694x_3^2x_5m_1m_4^2{-}12093x_1x_4x_5m_1m_4^2{-}12093x_2x_4x_5m_1m_4^2{-}12093x_3x_4x_5m_1m_4^2\\
{-}6694x_4^2x_5m_1m_4^2{-}6694x_1x_5^2m_1m_4^2{-}6694x_2x_5^2m_1m_4^2{-}6694x_3x_5^2m_1m_4^2\\
{-}6694x_4x_5^2m_1m_4^2{-}2663x_5^3m_1m_4^2{+}1748x_1^2m_1^2m_4^2{+}3190x_1x_2m_1^2m_4^2\\
{+}1748x_2^2m_1^2m_4^2{+}3190x_1x_3m_1^2m_4^2{+}3190x_2x_3m_1^2m_4^2{+}1748x_3^2m_1^2m_4^2\\
{+}3190x_1x_4m_1^2m_4^2{+}3190x_2x_4m_1^2m_4^2{+}3190x_3x_4m_1^2m_4^2{+}1748x_4^2m_1^2m_4^2\\
{+}3190x_1x_5m_1^2m_4^2{+}3190x_2x_5m_1^2m_4^2{+}3190x_3x_5m_1^2m_4^2{+}3190x_4x_5m_1^2m_4^2\\
{+}1748x_5^2m_1^2m_4^2{-}402x_1m_1^3m_4^2{-}402x_2m_1^3m_4^2{-}402x_3m_1^3m_4^2{-}402x_4m_1^3m_4^2\\
{-}402x_5m_1^3m_4^2{+}21m_1^4m_4^2{-}2663x_1^3m_2m_4^2{-}6694x_1^2x_2m_2m_4^2{-}6694x_1x_2^2m_2m_4^2\\
{-}2663x_2^3m_2m_4^2{-}6694x_1^2x_3m_2m_4^2{-}12093x_1x_2x_3m_2m_4^2{-}6694x_2^2x_3m_2m_4^2\\
{-}6694x_1x_3^2m_2m_4^2{-}6694x_2x_3^2m_2m_4^2{-}2663x_3^3m_2m_4^2{-}6694x_1^2x_4m_2m_4^2\\
{-}12093x_1x_2x_4m_2m_4^2{-}6694x_2^2x_4m_2m_4^2{-}12093x_1x_3x_4m_2m_4^2{-}12093x_2x_3x_4m_2m_4^2\\
{-}6694x_3^2x_4m_2m_4^2{-}6694x_1x_4^2m_2m_4^2{-}6694x_2x_4^2m_2m_4^2{-}6694x_3x_4^2m_2m_4^2\\
{-}2663x_4^3m_2m_4^2{-}6694x_1^2x_5m_2m_4^2{-}12093x_1x_2x_5m_2m_4^2{-}6694x_2^2x_5m_2m_4^2\\
{-}12093x_1x_3x_5m_2m_4^2{-}12093x_2x_3x_5m_2m_4^2{-}6694x_3^2x_5m_2m_4^2{-}12093x_1x_4x_5m_2m_4^2\\
{-}12093x_2x_4x_5m_2m_4^2{-}12093x_3x_4x_5m_2m_4^2{-}6694x_4^2x_5m_2m_4^2{-}6694x_1x_5^2m_2m_4^2\\
{-}6694x_2x_5^2m_2m_4^2{-}6694x_3x_5^2m_2m_4^2{-}6694x_4x_5^2m_2m_4^2{-}2663x_5^3m_2m_4^2\\
{+}4452x_1^2m_1m_2m_4^2{+}8122x_1x_2m_1m_2m_4^2{+}4452x_2^2m_1m_2m_4^2{+}8122x_1x_3m_1m_2m_4^2\\
{+}8122x_2x_3m_1m_2m_4^2{+}4452x_3^2m_1m_2m_4^2{+}8122x_1x_4m_1m_2m_4^2{+}8122x_2x_4m_1m_2m_4^2\\
{+}8122x_3x_4m_1m_2m_4^2{+}4452x_4^2m_1m_2m_4^2{+}8122x_1x_5m_1m_2m_4^2{+}8122x_2x_5m_1m_2m_4^2\\
{+}8122x_3x_5m_1m_2m_4^2{+}8122x_4x_5m_1m_2m_4^2{+}4452x_5^2m_1m_2m_4^2{-}2158x_1m_1^2m_2m_4^2\\
{-}2158x_2m_1^2m_2m_4^2{-}2158x_3m_1^2m_2m_4^2{-}2158x_4m_1^2m_2m_4^2{-}2158x_5m_1^2m_2m_4^2\\
{+}279m_1^3m_2m_4^2{+}1748x_1^2m_2^2m_4^2{+}3190x_1x_2m_2^2m_4^2{+}1748x_2^2m_2^2m_4^2\\
{+}3190x_1x_3m_2^2m_4^2{+}3190x_2x_3m_2^2m_4^2{+}1748x_3^2m_2^2m_4^2{+}3190x_1x_4m_2^2m_4^2\\
{+}3190x_2x_4m_2^2m_4^2{+}3190x_3x_4m_2^2m_4^2{+}1748x_4^2m_2^2m_4^2{+}3190x_1x_5m_2^2m_4^2\\
{+}3190x_2x_5m_2^2m_4^2{+}3190x_3x_5m_2^2m_4^2{+}3190x_4x_5m_2^2m_4^2{+}1748x_5^2m_2^2m_4^2\\
{-}2158x_1m_1m_2^2m_4^2{-}2158x_2m_1m_2^2m_4^2{-}2158x_3m_1m_2^2m_4^2{-}2158x_4m_1m_2^2m_4^2\\
{-}2158x_5m_1m_2^2m_4^2{+}593m_1^2m_2^2m_4^2{-}402x_1m_2^3m_4^2{-}402x_2m_2^3m_4^2{-}402x_3m_2^3m_4^2\\
{-}402x_4m_2^3m_4^2{-}402x_5m_2^3m_4^2{+}279m_1m_2^3m_4^2{+}21m_2^4m_4^2{-}2663x_1^3m_3m_4^2\\
{-}6694x_1^2x_2m_3m_4^2{-}6694x_1x_2^2m_3m_4^2{-}2663x_2^3m_3m_4^2{-}6694x_1^2x_3m_3m_4^2\\
{-}12093x_1x_2x_3m_3m_4^2{-}6694x_2^2x_3m_3m_4^2{-}6694x_1x_3^2m_3m_4^2{-}6694x_2x_3^2m_3m_4^2\\
{-}2663x_3^3m_3m_4^2{-}6694x_1^2x_4m_3m_4^2{-}12093x_1x_2x_4m_3m_4^2{-}6694x_2^2x_4m_3m_4^2\\
{-}12093x_1x_3x_4m_3m_4^2{-}12093x_2x_3x_4m_3m_4^2{-}6694x_3^2x_4m_3m_4^2{-}6694x_1x_4^2m_3m_4^2\\
{-}6694x_2x_4^2m_3m_4^2{-}6694x_3x_4^2m_3m_4^2{-}2663x_4^3m_3m_4^2{-}6694x_1^2x_5m_3m_4^2\\
{-}12093x_1x_2x_5m_3m_4^2{-}6694x_2^2x_5m_3m_4^2\\
{-}12093x_1x_3x_5m_3m_4^2{-}12093x_2x_3x_5m_3m_4^2{-}6694x_3^2x_5m_3m_4^2\\
{-}12093x_1x_4x_5m_3m_4^2{-}12093x_2x_4x_5m_3m_4^2\\
{-}12093x_3x_4x_5m_3m_4^2{-}6694x_4^2x_5m_3m_4^2{-}6694x_1x_5^2m_3m_4^2{-}6694x_2x_5^2m_3m_4^2{-}6694x_3x_5^2m_3m_4^2\\
{-}6694x_4x_5^2m_3m_4^2{-}2663x_5^3m_3m_4^2{+}4452x_1^2m_1m_3m_4^2{+}8122x_1x_2m_1m_3m_4^2{+}4452x_2^2m_1m_3m_4^2\\
{+}8122x_1x_3m_1m_3m_4^2{+}8122x_2x_3m_1m_3m_4^2{+}4452x_3^2m_1m_3m_4^2\\
{+}8122x_1x_4m_1m_3m_4^2{+}8122x_2x_4m_1m_3m_4^2\\
{+}8122x_3x_4m_1m_3m_4^2{+}4452x_4^2m_1m_3m_4^2{+}8122x_1x_5m_1m_3m_4^2\\
{+}8122x_2x_5m_1m_3m_4^2{+}8122x_3x_5m_1m_3m_4^2\\
{+}8122x_4x_5m_1m_3m_4^2{+}4452x_5^2m_1m_3m_4^2{-}2158x_1m_1^2m_3m_4^2{-}2158x_2m_1^2m_3m_4^2{-}2158x_3m_1^2m_3m_4^2\\
{-}2158x_4m_1^2m_3m_4^2{-}2158x_5m_1^2m_3m_4^2{+}279m_1^3m_3m_4^2{+}4452x_1^2m_2m_3m_4^2{+}8122x_1x_2m_2m_3m_4^2\\
{+}4452x_2^2m_2m_3m_4^2{+}8122x_1x_3m_2m_3m_4^2{+}8122x_2x_3m_2m_3m_4^2\\
{+}4452x_3^2m_2m_3m_4^2{+}8122x_1x_4m_2m_3m_4^2{+}8122x_2x_4m_2m_3m_4^2\\
{+}8122x_3x_4m_2m_3m_4^2{+}4452x_4^2m_2m_3m_4^2\\
{+}8122x_1x_5m_2m_3m_4^2{+}8122x_2x_5m_2m_3m_4^2{+}8122x_3x_5m_2m_3m_4^2\\
{+}8122x_4x_5m_2m_3m_4^2{+}4452x_5^2m_2m_3m_4^2\\
{-}5394x_1m_1m_2m_3m_4^2{-}5394x_2m_1m_2m_3m_4^2{-}5394x_3m_1m_2m_3m_4^2{-}5394x_4m_1m_2m_3m_4^2\\
{-}5394x_5m_1m_2m_3m_4^2{+}1456m_1^2m_2m_3m_4^2{-}2158x_1m_2^2m_3m_4^2{-}2158x_2m_2^2m_3m_4^2\\
{-}2158x_3m_2^2m_3m_4^2{-}2158x_4m_2^2m_3m_4^2\\
{-}2158x_5m_2^2m_3m_4^2{+}1456m_1m_2^2m_3m_4^2{+}279m_2^3m_3m_4^2{+}1748x_1^2m_3^2m_4^2{+}3190x_1x_2m_3^2m_4^2\\
{+}1748x_2^2m_3^2m_4^2{+}3190x_1x_3m_3^2m_4^2{+}3190x_2x_3m_3^2m_4^2{+}1748x_3^2m_3^2m_4^2{+}3190x_1x_4m_3^2m_4^2\\
{+}3190x_2x_4m_3^2m_4^2{+}3190x_3x_4m_3^2m_4^2{+}1748x_4^2m_3^2m_4^2{+}3190x_1x_5m_3^2m_4^2{+}3190x_2x_5m_3^2m_4^2\\
{+}3190x_3x_5m_3^2m_4^2{+}3190x_4x_5m_3^2m_4^2{+}1748x_5^2m_3^2m_4^2{-}2158x_1m_1m_3^2m_4^2{-}2158x_2m_1m_3^2m_4^2\\
{-}2158x_3m_1m_3^2m_4^2{-}2158x_4m_1m_3^2m_4^2{-}2158x_5m_1m_3^2m_4^2{+}593m_1^2m_3^2m_4^2{-}2158x_1m_2m_3^2m_4^2\\
{-}2158x_2m_2m_3^2m_4^2{-}2158x_3m_2m_3^2m_4^2{-}2158x_4m_2m_3^2m_4^2{-}2158x_5m_2m_3^2m_4^2{+}1456m_1m_2m_3^2m_4^2\\
{+}593m_2^2m_3^2m_4^2{-}402x_1m_3^3m_4^2{-}402x_2m_3^3m_4^2{-}402x_3m_3^3m_4^2{-}402x_4m_3^3m_4^2\\
{-}402x_5m_3^3m_4^2{+}279m_1m_3^3m_4^2{+}279m_2m_3^3m_4^2{+}21m_3^4m_4^2{-}459x_1^3m_4^3\\
{-}1152x_1^2x_2m_4^3{-}1152x_1x_2^2m_4^3{-}459x_2^3m_4^3{-}1152x_1^2x_3m_4^3\\
{-}2079x_1x_2x_3m_4^3{-}1152x_2^2x_3m_4^3{-}1152x_1x_3^2m_4^3{-}1152x_2x_3^2m_4^3{-}459x_3^3m_4^3\\
{-}1152x_1^2x_4m_4^3{-}2079x_1x_2x_4m_4^3{-}1152x_2^2x_4m_4^3{-}2079x_1x_3x_4m_4^3{-}2079x_2x_3x_4m_4^3\\
{-}1152x_3^2x_4m_4^3{-}1152x_1x_4^2m_4^3{-}1152x_2x_4^2m_4^3{-}1152x_3x_4^2m_4^3{-}459x_4^3m_4^3\\
{-}1152x_1^2x_5m_4^3{-}2079x_1x_2x_5m_4^3{-}1152x_2^2x_5m_4^3{-}2079x_1x_3x_5m_4^3{-}2079x_2x_3x_5m_4^3\\
{-}1152x_3^2x_5m_4^3{-}2079x_1x_4x_5m_4^3{-}2079x_2x_4x_5m_4^3{-}2079x_3x_4x_5m_4^3{-}1152x_4^2x_5m_4^3{-}1152x_1x_5^2m_4^3\\
{-}1152x_2x_5^2m_4^3{-}1152x_3x_5^2m_4^3{-}1152x_4x_5^2m_4^3{-}459x_5^3m_4^3{+}798x_1^2m_1m_4^3{+}1458x_1x_2m_1m_4^3\\
{+}798x_2^2m_1m_4^3{+}1458x_1x_3m_1m_4^3{+}1458x_2x_3m_1m_4^3{+}798x_3^2m_1m_4^3{+}1458x_1x_4m_1m_4^3\\
{+}1458x_2x_4m_1m_4^3{+}1458x_3x_4m_1m_4^3{+}798x_4^2m_1m_4^3{+}1458x_1x_5m_1m_4^3{+}1458x_2x_5m_1m_4^3\\
{+}1458x_3x_5m_1m_4^3{+}1458x_4x_5m_1m_4^3{+}798x_5^2m_1m_4^3{-}402x_1m_1^2m_4^3{-}402x_2m_1^2m_4^3{-}402x_3m_1^2m_4^3\\
{-}402x_4m_1^2m_4^3{-}402x_5m_1^2m_4^3{+}54m_1^3m_4^3{+}798x_1^2m_2m_4^3{+}1458x_1x_2m_2m_4^3{+}798x_2^2m_2m_4^3\\
{+}1458x_1x_3m_2m_4^3{+}1458x_2x_3m_2m_4^3{+}798x_3^2m_2m_4^3{+}1458x_1x_4m_2m_4^3{+}1458x_2x_4m_2m_4^3\\
{+}1458x_3x_4m_2m_4^3{+}798x_4^2m_2m_4^3{+}1458x_1x_5m_2m_4^3{+}1458x_2x_5m_2m_4^3{+}1458x_3x_5m_2m_4^3\\
{+}1458x_4x_5m_2m_4^3{+}798x_5^2m_2m_4^3{-}1002x_1m_1m_2m_4^3{-}1002x_2m_1m_2m_4^3{-}1002x_3m_1m_2m_4^3\\
{-}1002x_4m_1m_2m_4^3{-}1002x_5m_1m_2m_4^3{+}279m_1^2m_2m_4^3{-}402x_1m_2^2m_4^3\\
{-}402x_2m_2^2m_4^3{-}402x_3m_2^2m_4^3{-}402x_4m_2^2m_4^3{-}402x_5m_2^2m_4^3{+}279m_1m_2^2m_4^3\\
{+}54m_2^3m_4^3{+}798x_1^2m_3m_4^3{+}1458x_1x_2m_3m_4^3{+}798x_2^2m_3m_4^3{+}1458x_1x_3m_3m_4^3\\
{+}1458x_2x_3m_3m_4^3{+}798x_3^2m_3m_4^3{+}1458x_1x_4m_3m_4^3{+}1458x_2x_4m_3m_4^3{+}1458x_3x_4m_3m_4^3\\
{+}798x_4^2m_3m_4^3{+}1458x_1x_5m_3m_4^3{+}1458x_2x_5m_3m_4^3{+}1458x_3x_5m_3m_4^3{+}1458x_4x_5m_3m_4^3\\
{+}798x_5^2m_3m_4^3{-}1002x_1m_1m_3m_4^3{-}1002x_2m_1m_3m_4^3{-}1002x_3m_1m_3m_4^3\\
{-}1002x_4m_1m_3m_4^3{-}1002x_5m_1m_3m_4^3{+}279m_1^2m_3m_4^3{-}1002x_1m_2m_3m_4^3\\
{-}1002x_2m_2m_3m_4^3{-}1002x_3m_2m_3m_4^3{-}1002x_4m_2m_3m_4^3{-}1002x_5m_2m_3m_4^3\\
{+}684m_1m_2m_3m_4^3{+}279m_2^2m_3m_4^3{-}402x_1m_3^2m_4^3{-}402x_2m_3^2m_4^3{-}402x_3m_3^2m_4^3\\
{-}402x_4m_3^2m_4^3{-}402x_5m_3^2m_4^3{+}279m_1m_3^2m_4^3{+}279m_2m_3^2m_4^3{+}54m_3^3m_4^3{+}54x_1^2m_4^4\\
{+}99x_1x_2m_4^4{+}54x_2^2m_4^4{+}99x_1x_3m_4^4{+}99x_2x_3m_4^4{+}54x_3^2m_4^4{+}99x_1x_4m_4^4{+}99x_2x_4m_4^4\\
{+}99x_3x_4m_4^4{+}54x_4^2m_4^4{+}99x_1x_5m_4^4{+}99x_2x_5m_4^4{+}99x_3x_5m_4^4{+}99x_4x_5m_4^4\\
{+}54x_5^2m_4^4{-}72x_1m_1m_4^4{-}72x_2m_1m_4^4{-}72x_3m_1m_4^4{-}72x_4m_1m_4^4\\
{-}72x_5m_1m_4^4{+}21m_1^2m_4^4{-}72x_1m_2m_4^4{-}72x_2m_2m_4^4{-}72x_3m_2m_4^4\\
{-}72x_4m_2m_4^4{-}72x_5m_2m_4^4{+}51m_1m_2m_4^4{+}21m_2^2m_4^4{-}72x_1m_3m_4^4\\
{-}72x_2m_3m_4^4{-}72x_3m_3m_4^4{-}72x_4m_3m_4^4{-}72x_5m_3m_4^4{+}51m_1m_3m_4^4\\
{+}51m_2m_3m_4^4{+}21m_3^2m_4^4)
$
\end{document}